\documentclass[11pt,titlepage,a4paper]{article}
\usepackage{bozhomac}   
\usepackage{bozhlogo}   
\usepackage{cite}   
\usepackage{tabularx}   

\title{\bfseries    \vspace*{-1.678902345in}
{\huge Normal frames for general connections
\\[1.22ex]  on differentiable fibre bundles         }
}

\vspace{1.7ex}

\author{
Bozhidar Z.\ Iliev
\thanks{Laboratory of Mathematical Modeling in Physics,
Institute for Nuclear Research and \mbox{Nuclear} Energy,
Bulgarian Academy of Sciences,
Boul.\ Tzarigradsko chauss\'ee~72, 1784 Sofia, Bulgaria}
\thanks{E-mail address: bozho@inrne.bas.bg}
\thanks{URL: http://theo.inrne.bas.bg/$\sim$bozho/}
}

\date{
 \vspace{2.27ex}\ShortTitle{Normal frames for connections}  \\[0.27ex]
 \vspace{3.27ex}
\small
    \begin{tabular}{r@{$\colon\to~$}l}
 \vspace{0.09ex} Began          & June 13, 2003     \\[0.09ex]
 \vspace{0.09ex} Ended          & August 16, 2003     \\[0.09ex]
\vspace{0.09ex} Initial typeset& July 22-- August 22, 2003 \\[0.09ex]
 \vspace{0.09ex} Last update   & May 31, 2004    \\[0.09ex]
 \vspace{0.27ex} Produced   & \fbox{\today} \\[0.27ex]
    \end{tabular} \\[1.27ex]
\normalsize
    \begin{tabular}{r@{$\colon~$}l}
 \vspace{0.27ex} http://www.arXiv.org e-Print archive No. & math.DG/0405004
                                \\[0.27ex]
    \end{tabular} \\[-0.27ex]
 \vspace{4.27ex}{\Huge \BOZHO}  \\[4.27ex]
    \begin{tabular}{r@{\hspace{0.512em}}|@{\hspace{0.512em}}l}
 \vspace{0.27ex}
&
 \vspace{0.27ex}
    \end{tabular} \\[1.27ex]
 \vspace{0.27ex}
					\\[0.27ex]
}

\begin{document}        

\renewcommand{\thepage}{\roman{page}}

\renewcommand{\thefootnote}{\fnsymbol{footnote}} 
\maketitle              
\renewcommand{\thefootnote}{\arabic{footnote}}   

\tableofcontents        


\begin{abstract}

The theory of frames normal for general connections on differentiable bundles
is developed. Links with the existing theory of frames normal for
covariant derivative operators (linear connections) in vector bundles are
revealed. The existence of bundle coordinates normal at a given point and/or
along injective horizontal path is proved. A necessary and sufficient
condition of existence of bundle coordinates normal along injective
horizontal mappings is derived.

\end{abstract}

\renewcommand{\thepage}{\arabic{page}}


\section {Introduction}
\label{Introduction}

    Generally said, normal are called frames and coordinates in which
the components of some geometrical object (locally) vanish; in particular,
these can be the coefficients of a linear connection on a manifold or in
vector bundle. As a result of that, some objects look like in a ``flat'' or
``Euclidean'' case, which significantly simplifies certain calculations,
formulae, their interpretation, etc. For instance, the normal frames for
linear connections turn to be the mathematical object for description of the
inertial frames of reference in physics, in which some effects of a force
field, like the gravity one, locally disappear.

    The history of the theory of normal coordinates and frames goes back
to 1854. The major classical results concerning the normal coordinates for
linear connections are summarize in the table below.

\renewcommand{\arraystretch}{1.32}
    \begin{table}[ht!]  \label{HistoricalTable}%
\index{normal coordinates!history of}
    \begin{tabularx}{\textwidth}{@{}rlX@{}}
Year    & Person    & Result and original reference
\\ \hline
1854    & B.~Riemann%
            & Existence and construction of (`Riemannian')
coordinates in a Riemannian manifold which are normal at a single
point.~\cite{Riemann/Hypotheses}
\\
1922    & O.~Veblen%
            & Existence and construction of (`Riemannian normal')
coordinates in a manifold with torsionless linear connection which are normal
at a single point.~\cite{Veblen}
\\
1922    & E.~Fermi%
            & Existence of (`Fermi') coordinates in a Riemannian
manifold which are normal along a path without self\ndash
intersections.~\cite{Fermi}
\\
1926    & T.~Levi-Civita%
            & Explicit transformation to the Fermi coordinates
along paths without self\ndash intersections.~\cite{Levi-Civita/1926}
\\
1927    & L.~P.~Eisenhart%
            & Existence and construction of particular
kind of (`Fermi') coordinates on a manifold with torsionless linear
connection which are normal along a path without self\ndash
intersections.~\cite{Eisenhart/Non-Riemannian}
\\
1958    & L.~O'Raifeartaigh%
            & Necessary and sufficient conditions for existence
of coordinates normal on submanifold of a manifold with torsionless linear
connection. If such coordinates exist, a particular example of them (`Fermi
coordinates') is constructed.~\cite{ORai}
\\\hline
    \end{tabularx}
    \end{table}

    In~\cite{bp-Frames-n+point,bp-Frames-path,bp-Frames-general} the
normal frames were introduced and studied for derivations, in particular for
linear connections, with generally non\ndash vanishing curvature and torsion
on a differentiable manifold. Then these objects were investigated for
derivations and linear connections in vector bundles~\cite{bp-NF-D+EP}. At
last, the paper~\cite{bp-NF-LTP} explores them for linear transports along
paths in vector bundles. The present work is devoted to the introduction and
some properties of normal frames and coordinates for general connections on
fibre bundles whose bundle and base spaces are differentiable manifolds.

    The lay-out of the work is as follows.

    In Sect.~\ref{Sect2} is collected some introductory material needed
for our exposition. Here some of our notation is fixed too.

    Section~\ref{Sect3} is devoted to the general connection theory on
bundles whose base and bundles spaces are differentiable manifolds. In
Subsect.~\ref{Subsect3.1} are reviewed some coordinates and frames/bases on
the bundle space which are compatible with the fibre structure of a bundle.
Subsect.~\ref{Subsect3.2} deals with the general connection theory. A connection
on a bundle is defined as a distribution on its bundle space which is
complimentary to the vertical distribution on it. The notion of specialized
frame is introduced. Frames adapted to  specialized frames, in
particular to local bundle coordinates, are defined and the local
(2\ndash index) coefficients in them of a connection are defined and their
transformation law is derived.

    The theory of normal frames for connections on bundles is considered
in sections~\ref{Sect7.1}--\ref{Sect7.3}. Sect.~\ref{Sect7.1} deals with
the general case. Loosely said, an adapted frame is called normal if in it
vanish the 2\ndash index coefficients of a connection on some set. It happens
that a frame is normal if and only if it coincides with the frame it is
adapted to. The set of these frames is completely described in the most
general case. The problems of existence, uniqueness, \etc of normal frames
adapted to holonomic frames, \ie adapted to local coordinates, are discussed
in Sect.~\ref{Sect7.2}. If such frames exist, their general form is
described. The existence of frames normal at a given point and/or
along an injective horizontal path is proved. The flatness of a connection
on an open set is pointed as a necessary condition of existence of (locally)
holonomic frames normal on that set. Some links between the general theory of
normal frames and the existing one of normal frames in vector bundles are given
in Sect.~\ref{Sect7.3}. It is proved that a frame is normal on a vector bundle
with linear connection if and only if in it vanish the 3\ndash index
coefficients of the connection. The equivalence of the both theories on vector
bundles is established.

    Section~\ref{Conclusion} ends the paper with some concluding
remarks.

	In the Appendix is formulated and proved a necessary and sufficient
condition for the existence of coordinates normal along injective mappings with
non\ndash vanishing horizontal component, in particular along injective
horizontal mappings.


\section {Preliminaries}
\label{Sect2}

        This section contains an introductory material, notation etc.\ that
will be needed for our exposition. The reader is referred for details to
standard books on differential geometry, like~\cite{K&N,Warner,Poor}.

        A differentiable finite-dimensional manifold over a field $\field$
will be denoted typically by $M$. Here $\field$ stands for the field
$\field[R]$ of real or the field $\field[C]$ of complex numbers,
$\field=\field[R],\field[C]$. The manifolds we consider are supposed to be
smooth of class $C^2$.~%
\footnote{~%
Some of our definitions or/and results are valid also for $C^1$ or even $C^0$
manifolds, but we do not want to overload the material with continuous
counting of the required degree of differentiability of the manifolds
involved. Some parts of the text admit generalizations on more general
spaces, like the topological ones, but this is out of the subject of the
present work.%
}
The sets of vector fields, realized as first order differential operators,
over $M$ will be denoted by $\mathcal{X}(M)$. The space tangent to $M$ at
$p\in M$ is $T_p(M)$ and $(T(M),\pi_T,M)$ will stand for the tangent bundle
over $M$. The value of $X\in\mathcal{X}(M)$ at $p\in M$ is $X_p\in T_p(M)$.

        If $M$ and $\bar{M}$ are manifolds and $f\colon \bar{M}\to M$ is  a
$C^1$ mapping, then $f_*:=\od f \colon T(\bar{M})\to T(M)$ denotes the
induced tangent mapping (or differential) of $f$ such that, for $p\in M$,
 $f_*|_p :=\od f|_p \colon T_p(\bar{M})\to T_{f(p)}(M)$
and, for a $C^1$ function $g$ on $M$,
 $(f_*(X))(g):=X(g\circ f)\colon p\mapsto f_*|_p(g)=X_p(g\circ f)$, with
$\circ$ being the composition of mappings sign.

    By $J\subseteq\field[R]$ will be denoted an arbitrary real interval
that can be open or closed at one or both its ends.
	The notation $\gamma\colon J\to M$ represents an arbitrary path in $M$.
	For a $C^1$ path $\gamma\colon J\to M$, the vector tangent to
$\gamma$ at $s\in J $ will be denoted by
\(
\dot\gamma(s) :=\frac{\od} {\od r}\Big|_{r=s}(\gamma(r))
=\gamma_*\bigl( \frac{\od}{\od r} \big|_{r=s} \bigr) \in T_{\gamma(s)}(M) .
\)
If $s_0\in J$ is an end point of $J$ and $J$ is
closed at $s_0$, the derivative in the definition of $\dot{\gamma}(s_0)$ is
regarded as a one\ndash sided derivative at $s_0$.

        Let the Greek indices $\lambda,\mu,\nu,\dots$ run over the range
$1,\dots,\dim M$ and $\{E_\mu\}$ be a $C^1$ frame in $T(M)$, \ie
$E_\mu\in\mathcal{X}(M)$ be of class $C^1$ and, for each $p\in M$, the set
$\{E_\mu|_p\}$ to be a basis of the vector space $T_p(M)$.~%
\footnote{~%
There are manifolds, like the even-dimensional spheres $\mathbb{S}^{2k}$,
$k\in\field[N]$, which do not admit global, continuous (and moreover $C^k$
for $k\ge1$), and nowhere vanishing vector fields~\cite{Spivak-1}. If this
is the case, the considerations must be localized over an open subset of $M$
on which such fields exist. We shall not overload our exposition with such
details.%
}
The Einstein's summation convention, summation on indices repeated on
different levels over the whole range of their values, will be assumed
hereafter.

        A frame $\{E_\mu\}$ or its dual coframe $\{E^\mu\}$ is called
\emph{holonomic} (\emph{anholonomic}) if $C_{\mu\nu}^{\lambda}=0$
($C_{\mu\nu}^{\lambda}\not=0$) for all (some) values of the indices $\mu$,
$\nu$, and $\lambda$, where the functions
$C_{\mu\lambda}^{\nu}$ are defined by
\(
[E_\mu,E_\nu]_{\_}
:= E_\mu\circ E_\nu - E_\nu\circ E_\mu
=: C_{\mu\nu}^{\lambda} E_\lambda ;
\)
these functions are a measure of deviation from a holonomic frame and are
known as the \emph{components of the anholonomy object of} $\{E_\mu\}$.
	For a holonomic frame always exist local coordinates
$\{x^\mu\}$ on $M$ such that \emph{locally} $E_\mu=\frac{\pd}{\pd x^\mu}$ and
$E^\mu=\od x^\mu$. Conversely, if $\{x^\mu\}$ are local coordinates on $M$,
then the local frame $\bigl\{\frac{\pd}{\pd x^\mu}\bigr\}$ and local coframe
$\{\od x^\mu\}$ are defined and holonomic on the domain of $\{x^\mu\}$.

        If $n\in \field[N]$ and $n\le\dim M$, an  $n$-dimensional
\emph{distribution} $\Delta$ on $M$ is defined as a mapping
$\Delta\colon p\mapsto \Delta_p$ assigning to each $p\in M$ an $n$\ndash
dimensional subspace $\Delta_p$ of the tangent space $T_p(M)$ of $M$ at $p$,
$\Delta_p\subseteq T_p(M)$.
    A distribution is \emph{integrable} if there is a submersion
$\psi\colon M\to N$ such that $\Ker \psi_*=\Delta$; a necessary and locally
sufficient condition for the integrability of $\Delta$ is the commutator of
every two vector fields in $\Delta$ to be in $\Delta$. We say that a vector
field $X\in\mathcal{X}(M)$ is in $\Delta$ and write $X\in\Delta$, if
$X_p\in\Delta_p$ for all $p\in M$. A \emph{basis on $U\subseteq M$ for}
$\Delta$ is a set  $\{X_1,\dots,X_n\}$ of $n$ linearly independent (relative to
functions $U\to\field$) vector fields in $\Delta|_U$, \ie
$\{X_1|_p,\dots,X_n|_p\}$ is a basis for $\Delta_p$ for all $p\in U$.

        A distribution is convenient to be described in terms of (global)
frames or/and coframes over $M$. In fact, if $p\in M$ and
$\varrho=1,\dots,n$, in each $\Delta_p\subseteq T_p(M)$, we can choose a basis
$\{X_\varrho|_p\}$ and hence a frame $\{X_\varrho\}$,
$X_\varrho\colon p\mapsto X_\varrho|_p$, in
$\{\Delta_p : p\in M\} \subseteq T(M) $; we say that $\{X_\rho\}$ is a basis
for/in $\Delta$. Conversely, any collection of $n$ linearly independent
(relative to functions $M\to\field$) vector fields $X_\varrho$ on $M$ defines
a distribution
\(
 p\mapsto\bigl\{ \sum_{\varrho=1}^{n} f^\varrho X_\varrho|_p :
f^\varrho\in\field \bigr\}.
\)
Consequently, a frame in $T(M)$ can be formed by
adding to a basis for $\Delta$ a set of $(\dim M-n)$ new linearly independent
vector fields (forming a frame in $T(M)\setminus\{\Delta_p : p\in M\}$)  and
v.v., by selecting  $n$ linearly independent vector fields on $M$, we can
define a distribution $\Delta$ on $M$.


\section {Connections on bundles}
\label{Sect3}

    Before presenting the general connection theory in
Subsect.~\ref{Subsect3.2}, we at first fix some notation and concepts
concerning fibre bundles in Subsect.~\ref{Subsect3.1}.


\subsection{Coordinates and frames on the bundle space}
\label{Subsect3.1}

         Let $(E,\pi,M$) be a bundle with bundle space $E$, projection
$\pi\colon E\to M$, and base space $M$. We suppose that the spaces $E$ and $M$
are $C^2$ differentiable, if the opposite is not stated explicitly,~%
\footnote{~%
Most of our considerations are valid also if $C^1$ differentiability is
assumed and even some of them hold on $C^0$ manifolds. By assuming $C^2$
differentiability, we skip the problem of counting the required
differentiability class of the whole material that follows. Sometimes, the
$C^2$  differentiability is required explicitly, which is a hint that a
statement or definition is not valid otherwise. If we want to emphasize that
some text is valid under a $C^1$ differentiability assumption, we indicate
that fact explicitly. However, the proofs of lemmas~\ref{2-Lem6.1-0}
and~\ref{3-Lem8.1-0}, proposition~\ref{Prop7.5} and all assertions in
section~\ref{Appendix} require $C^3$ differentiability, which will be indicated
explicitly%
}
manifolds of finite dimensions $n\in\field[N]$ and $n+r$, for some
$r\in\field[N]$, respectively; so the dimension of the fibres $\pi^{-1}(x)$,
with $x\in M$, \ie the fibre dimensions of $(E,\pi,M)$, is $r$.

         Let the Greek indices $\lambda,\mu,\nu,\ldots$ run from 1 to
$n=\dim M$, the Latin indices $a,b,c,\ldots$  take  the  values  from  $n+1$
to  $n+r=\dim E$, and the uppercase Latin indices
${I},{J},{K},\ldots$ take values in the whole set $\{1,\ldots,n+r\}$. One may
call these types of indices respectively base, fibre, and bundle indices.

        Suppose $\{u^{I}\}=\{u^\mu,u^a\}=\{u^1,\dots,u^{n+r}\}$ are local
bundle coordinates on an open set $U\subseteq E$, \ie on the set
$\pi(U)\subseteq M$ there are local coordinates $\{x^\mu\}$ such that
$u^\mu=x^\mu\circ \pi$; the coordinates $\{u^\mu\}$ (resp.\ $\{u^a\}$) are
called \emph{basic} (resp.\ \emph{fibre}) \emph{coordinates}~\cite{Poor}.~%
\footnote{~%
If $(U,v)$ is a bundle chart, with $v\colon U\to\field^n\times\field^r$ and
$e^a\colon\field^r\to\field$ are such that $e^a(c_1,\dots,c_r)=c_a\in\field$,
then one can put $u^a=e^a\circ \pr_2\circ v$, where
$\pr_2\colon\field^n\times\field^r\to\field^r$ is the projection on the second
multiplier $\field^r$.%
}

        Further only coordinate changes
        \begin{subequations}    \label{3.1}
        \begin{equation}        \label{3.1a}
\{u^\mu,u^a\} \mapsto \{\tilde{u}^\mu,\tilde{u}^a\}
        \end{equation}
on $E$ which respect the fibre structure, \viz the division into basic and
fibre coordinates, will be considered; this means that
        \begin{equation}        \label{3.1b}
        \begin{split}
\tilde{u}^\mu(p) &= f^\mu(u^1(p),\dots,u^n(p))\\
\tilde{u}^a(p)   &= f^a(u^1(p),\dots,u^n(p),u^{n+1}(p),\dots,u^{n+r}(p))
        \end{split}
        \end{equation}
        \end{subequations}
for $p\in E$ and some functions $f^{I}$. The bundle coordinates
$\{u^\mu,u^a\}$ induce the (local) frame
$\bigl\{ \pd_ \mu:=\frac{\pd}{\pd u^\mu}, \pd_a:=\frac{\pd}{\pd u^a} \bigr\}$
 over $U$ in the tangent bundle space $T(E)$ of the tangent bundle over the
bundle space $E$. Since a change~\eref{3.1} of the coordinates on $E$ implies
\(
\pd_{I}\mapsto\tilde{\pd}_{I}
:=\frac{\pd}{\pd\tilde{u}^{I}}
=\frac{\pd u^{J}}{\pd\tilde{u}^{I}} \pd_{J} ,
\)
the transformation~\eref{3.1} leads to
        \begin{align}   \label{3.2a}
(\pd_\mu,\pd_a) &\mapsto
(\tilde{\pd}_\mu,\tilde{\pd}_a) = (\pd_\nu,\pd_b) \cdot A
        \end{align}
Here expressions like $(\pd_\mu,\pd_a)$ are shortcuts for ordered
$(n+r)$\ndash tuples like
$(\pd_1,\dots,\pd_{n+r}) = \bigl([\pd_\mu]_{\mu=1}^n,[\pd_a]_{a=n+1}^{n+r}$),
the centered dot $\cdot$ stands
for the matrix multiplication, and the transformation matrix $A$ is
        \begin{equation}        \label{3.3}
A
:= \Bigl[\frac{\pd u^{I}}{\pd\tilde{u}^{J}}\Bigr]_{{I},{J}=1}^{n+r}
=
        \begin{pmatrix}
\bigl[\frac{\pd u^\nu}{\pd\tilde{u}^\mu}\bigr] & 0_{n\times r}
\\
\bigl[\frac{\pd u^b}{\pd\tilde{u}^\mu}\bigr] &
\bigl[\frac{\pd u^b}{\pd\tilde{u}^a}\bigr]
        \end{pmatrix}
=:
        \begin{bmatrix}
\frac{\pd u^\nu}{\pd\tilde{u}^\mu} & 0
\\
\frac{\pd u^b}{\pd\tilde{u}^\mu} &
\frac{\pd u^b}{\pd\tilde{u}^a}
        \end{bmatrix}
\ ,
        \end{equation}
where $0_{n\times r}$ is the $n\times r$ zero matrix. The explicit form of
the matrix inverse to~\eref{3.3} is
 $A^{-1}=\bigl[\frac{\pd\tilde{u}^{I}}{\pd u^{J}}\bigr]=\ldots$ and it
is obtained from~\eref{3.3} via the change $u\leftrightarrow\tilde{u}$.

        The formula~\eref{3.2a} can be generalized for arbitrary frame
$\{e_{I}\}=\{e_\mu,e_a\}$ in $T(E)$
which respect the fibre structure in
a sense that their \emph{admissible changes} are given by
        \begin{align}   \label{3.4a}
(e_{I})=(e_\mu,e_a)
&\mapsto
(\tilde{e}_{I})=(\tilde{e}_\mu,\tilde{e}_a)
= (e_\nu,e_b)\cdot A
        \end{align}
Here $A=[A_{J}^{I}]$ is a nondegenerate matrix-valued function with a block
structure similar to~\eref{3.3}, \viz
        \begin{subequations}        \label{3.5}
        \begin{align}   \label{3.5a}
A =
        \begin{pmatrix}
[A_\mu^\nu]_{\mu,\nu=1}^{n} & 0_{n\times r}
\\
\bigl[A_\mu^b\bigr]_{
           \begin{subarray}{l}
           \mu=1,\dots,n \\
           b=n+1,\dots,n+r
           \end{subarray}
          }
&
[A_a^b]_{a,b=n+1}^{n+r}
        \end{pmatrix}
=:
        \begin{bmatrix}
A_\mu^\nu & 0
\\
A_\mu^b   & A_a^b
        \end{bmatrix}
\\\intertext{with inverse matrix}
                        \label{3.5b}
A^{-1}
=
        \begin{pmatrix}
[A_\mu^\nu]^{-1} & 0
\\
[A_b^a]^{-1} \cdot [A_\mu^a]\cdot [A_\mu^\nu]^{-1}
&
[A_b^a]^{-1}
        \end{pmatrix}\ ,
        \end{align}
        \end{subequations}
Here $A_\mu^a\colon U\to\field$ and $[A_\mu^\nu]$ and $[A_b^a]^{-1}$
are non-degenerate matrix\ndash valued functions on $U$ such that
$[A_\mu^\nu]$ is constant on the fibres of $E$, i.e., for $p\in E$,
$A_\mu^\nu(p)$ depends only on $\pi(p)\in M$, which is equivalent to any one
of the equations
$ A_\mu^\nu = B_\mu^\nu \circ \pi $
and
$\frac{\pd A_\mu^\nu }{\pd u^a} = 0 ,$
with $[B_\mu^\nu]$ being a nondegenerate matrix-valued function on
$\pi(U)\subseteq M$. Obviously,~\eref{3.2a} corresponds
to~\eref{3.4a} with
$e_I=\frac{\pd}{\pd{u}^I}$,
$\tilde{e}_I=\frac{\pd}{\pd\tilde{u}^I}$, and
$A_{I}^{J}=\frac{\pd u^{J}}{\pd\tilde{u}^{I}}$.

	All frames on $E$ connected via~\eref{3.4a}--\eref{3.5}, which are
(locally) obtainable from holonomic ones, induced by bundle coordinates, via
admissible changes, will be referred as \emph{bundle frames}.



\subsection{Connection theory}
\label{Subsect3.2}

	From a number of equivalent definitions of a connection on
differentiable manifold~\cite[sections~2.1
and~2.2]{Mangiarotti&Sardanashvily}, we shall use the following one.

    \begin{Defn}    \label{Defn3.1}
A \emph{connection on a bundle} $(E,\pi,M)$ is an $n=\dim M$
dimensional distribution $\Delta^h$ on $E$ such that, for each $p\in E$
and the \emph{vertical distribution} $\Delta^v$ defined by
    \begin{equation}    \label{3.9-2}
\Delta^v\colon p \mapsto \Delta^v_p
:= T_{\imath(p)}\bigl( \pi^{-1}(\pi(p)) \bigr)
\cong T_{p}\bigl( \pi^{-1}(\pi(p)) \bigr),
    \end{equation}
with $\imath\colon\pi^{-1}(\pi(p))\to E$ being the inclusion mapping, is
fulfilled
    \begin{equation}    \label{3.9-3}
\Delta^v_p\oplus \Delta^h_p = T_p(E) ,
    \end{equation}
where
\(
\Delta^h\colon p \mapsto \Delta^h_p \subseteq
T_{p}\bigl( \pi^{-1}(\pi(p)) \bigr)
\)
and $\oplus$ is the direct sum sign. The distribution $\Delta^h$ is
called \emph{horizontal} and symbolically we write
$\Delta^v\oplus\Delta^h=T(E)$.
    \end{Defn}

    A \emph{vector} at a point $p\in E$ (resp. a \emph{vector field} on $E$)
is said to be \emph{vertical} or \emph{horizontal} if it (resp.\ its value at
$p$) belongs to $\Delta^h_p$ or $\Delta^v_p$, respectively, for the given
(resp.\ any) point $p$. A vector $Y_p\in T_p(E)$ (resp.\ vector field
$Y\in\mathcal{X}(E)$) is called a \emph{horizontal lift of a vector}
$X_{\pi(p)}\in T_{\pi(p)}(M)$ (resp.\ \emph{vector field} $X\in\mathcal{X}(M)$
on $M=\pi(E)$) if $\pi_*(Y_p)=X_{\pi(p)}$ for the given (resp.\ any) point
$p\in E$. Since $\pi_*|_{\Delta_p^h}\colon \Delta_p^h\to T_{\pi(p)}(M)$ is a
vector space isomorphism for all $p\in E$~\cite[sec.~1.24]{Poor}, any vector in
$T_{\pi(p)}(M)$ (resp.\ vector field in $\mathcal{X}(M)$) has a unique
horizontal lift in $T_p(E)$ (resp.\ $\mathcal{X}(E)$).

    As a result of~\eref{3.9-3}, any vector $Y_p\in T_p(E)$ (resp.\ vector
field $Y\in\mathcal{X}(E)$) admits a unique representation
 $Y_p=Y_p^v\oplus Y_p^h$ (resp.\ $Y=Y^v\oplus Y^h$) with
 $Y_p^v\in\Delta_p^v$ and $Y_p^h\in\Delta_p^h$
(resp.\ $Y^v\in\Delta^v$ and $Y^h\in\Delta^h$). If the distribution
$p\mapsto\Delta_p^h$ is differentiable of class $C^m$,
$m\in\field[N]\cup\{0,\infty,\omega\}$, it is said that the \emph{connection
$\Delta^h$ is (differentiable) of class} $C^m$. A connection $\Delta^h$ is of
class $C^m$ if and only if, for every $C^m$ vector field $Y$ on $E$, the
vertical $Y^v$ and horizontal $Y^h$ vector fields are of class $C^m$.

    Let us now look on a connections $\Delta^h$ on a bundle $(E,\pi,M)$
from a view point of (local) frames and their dual coframes on $E$. Let
$\{e_\mu\}$ be a basis for $\Delta^h$, \ie $e_\mu\in\Delta^h$ and
$\{e_\mu|_p\}$ is a basis for $\Delta_p^h$ for all $p\in E$.

    \begin{Defn}    \label{Defn3.3}
A frame $\{e_I\}$ in $T(E)$ over $E$ is called \emph{specialized} for a
connection $\Delta^h$ if the first $n=\dim M$ of its vector fields
$\{e_\mu\}$ form a basis for the horizontal distribution $\Delta^h$ and its
last $r=\dim\pi^{-1}(x)$, $x\in M$, vector fields $\{e_a\}$ form a basis for
the vertical distribution $\Delta^v$.
    \end{Defn}


    It is a simple, but important, fact that the specialized frames are
the most general ones which respect the splitting of $T(E)$ into vertical and
horizontal components. Suppose $\{e_{I}\}$ is a specialized frame. Then
the general element of the set of all specialized frames is (see~\eref{3.4a})
    \begin{equation}    \label{6.3a}
(\bar{e}_\mu,\bar{e}_a)
=
({e}_\nu,{e}_b) \cdot
    \begin{bmatrix}
A_\mu^\nu   & 0 \\
0       & A_a^b
    \end{bmatrix}
=
(A_\mu^\nu e_\nu, A_a^b e_b),
    \end{equation}
where $[A_\mu^\nu]_{\mu,\nu=1}^{n}$ and
$[A_a^b]_{a,b=n+1}^{n+r}$ are non-degenerate matrix-valued functions on $E$,
which are constant on the fibres of $(E,\pi,M)$, \ie we can set
 $A_\mu^\nu=B_\mu^\nu\circ\pi$ and $A_a^b=B_a^b\circ\pi$ for some
non\ndash degenerate matrix\ndash valued functions $[B_\mu^\nu]$ and
$[B_a^b]$ on $M$.

    Since $\pi_*|_{\Delta^h}\colon\{X\in\Delta^h\}\to \mathcal{X}(M)$ is
an isomorphism, any basis $\{\varepsilon_\mu\}$ \emph{for} $\Delta^h$
defines a basis $\{E_\mu\}$ of $\mathcal{X}(M)$ such that
    \begin{equation}    \label{6.4}
E_\mu = \pi_*|_{\Delta^h}(\varepsilon_\mu)
    \end{equation}
and v.v., a basis $\{E_\mu\}$ for $\mathcal{X}(M)$ induces a basis
$\{\varepsilon_\mu\}$ for $\Delta^h$ via
    \begin{equation}    \label{6.5}
\varepsilon_\mu = (\pi_*|_{\Delta^h})^{-1} (E_\mu) .
    \end{equation}
Thus a `horizontal' change
    \begin{equation}    \label{6.6}
\varepsilon_\mu\mapsto
\bar{\varepsilon}_\mu = (B_\mu^\nu\circ\pi) \varepsilon_\nu ,
    \end{equation}
which is independent of a `vertical' one given by
    \begin{equation}    \label{6.7}
\varepsilon_a\mapsto
\bar{\varepsilon}_a = (B_a^b\circ\pi) \varepsilon_b
    \end{equation}
with $\{\varepsilon_a\}$ being a basis for $\Delta^v$, is
equivalent to the transformation
    \begin{equation}    \label{6.8}
E_\mu \mapsto \bar{E}_\mu = B_\mu^\nu E_\nu
    \end{equation}
of the basis $\{E_\mu\}$ for $\mathcal{X}(M)$, related via~\eref{6.4} to the
basis $\{\varepsilon_\mu\}$ for $\Delta^h$. Here $[B_\mu^\nu]$ and $[B_a^b]$
are non\ndash degenerate matrix\ndash valued functions on $M$.

    As $\pi_*(\varepsilon_a)=0\in\mathcal{X}(M)$, the `vertical'
transformations~\eref{6.7} do not admit interpretation analogous to the
`horizontal' ones~\eref{6.6}. However, in a case of a \emph{vector} bundle
$(E,\pi,M)$, they are tantamount to changes of frames in the bundle space
$E$, \ie of the bases for $\Sec(E,\pi,M)$. Indeed, if $v$ is a mapping
defined by
    \begin{equation}    \label{4.1-1}
    \begin{split}
v & \colon\Sec(E,\pi,M)  \to \{ \text{vector fields in }\Delta^v \}
\\
v & \colon Y  \mapsto Y^v \colon p\mapsto
            Y^v|_p  := \frac{\od}{\od t}\Big|_{t=0}(p+tY_{\pi(p)}) ,
    \end{split}
    \end{equation}
which mapping is a linear isomorphism~\cite{Poor
},
the sections
    \begin{equation}    \label{6.9}
E_a = v^{-1}(\varepsilon_a)
    \end{equation}
form a basis for $\Sec(E,\pi,M)$ as the vertical vector fields
$\varepsilon_a$ form a basis for $\Delta^v$. Conversely, any basis $\{E_a\}$
for the sections of $(E,\pi,M)$ induces a basis $\{\varepsilon_a\}$ for
$\Delta^v$ such that
    \begin{equation}    \label{6.10}
\varepsilon_a = v(E_a) .
    \end{equation}
As $v$ and $v^{-1}$ are linear, the change~\eref{6.7} is equivalent to the
transformation
    \begin{equation}    \label{6.11}
E_a\mapsto \bar{E}_a = B_a^b E_b
    \end{equation}
of the frame $\{E_a\}$ in $E$ related to $\{\varepsilon_a\}$ via~\eref{6.9}.
In this way, we see that
\emph{any specialized frame
$\{\varepsilon_{I}\}=\{\varepsilon_\mu,\varepsilon_a\}$
for a connection on a vector bundle $(E,\pi,M)$ is equivalent to a pair of
frames $(\{E_\mu\},\{E_a\})$%
}
such that $\{E_\mu\}$ is a basis for the set $\mathcal{X}(M)$ of vector
fields on the base $M$, \ie for the sections of the tangent bundle
$(T(M),\pi_T,M)$ (and hence is a frame in $T(M)$ over $M$), and $\{E_a\}$ is
a basis for the set $\Sec(E,\pi,M)$ of sections of the initial bundle (and
hence is a frame in $E$ over $M$). Since conceptually the frames in $T(M)$ and
$E$ are easier to be understood and in some cases have a direct physical
interpretation, one often works with the pair of frames
\(
( \{E_\mu=\pi_*|_{\Delta^h}(\varepsilon_\mu)\} ,
\{E_a=v^{-1}(\varepsilon_a)\} )
\)
instead with a specialized frame
$\{\varepsilon_{I}\}=\{\varepsilon_\mu,\varepsilon_a\}$;
for instance $\{E_\mu\}$ and $\{E_a\}$ can be completely arbitrary frames in
$T(M)$ and $E$, respectively, while the specialized frames represent only a
particular class  of frames in $T(E)$.

    One can \emph{mutatis mutandis} localize the above considerations
when $M$ is replaced with an open subset $U_M$ in $M$ and $E$ is replaced with
$U=\pi^{-1}(U_M)$. Such a localization is important when the bases/frames
considered are connected with some local coordinates or when they should be
smooth.%
\footnote{~%
Recall, not every manifold admits a \emph{global} nowhere vanishing $C^m$,
$m\ge0$, vector field (see~\cite{Spivak-1} or~\cite[sec.~4.24]{Schutz}); \eg
such are the even\ndash dimensional spheres $\mathbb{S}^{2k}$,
$k\in\field[N]$, in Euclidean space.%
}

	Let $\{e_{I}\}$ be a frame in $T(E)$ defined over an open set
$U\subseteq E$ and such that $\{e_a|_p\}$ is a \emph{basis for the space
$T_p(\pi^{-1}(\pi(p)))$ tangent to the fibre through} $p\in U$. Then we can
write the expansion
    \begin{equation}    \label{6.12}
(e_\mu^U,e_a^U )
= (D_\mu^\nu e_\mu + D_\mu^a e_a , D_a^b e_b )
=
(e_\nu,e_b) \cdot
    \begin{pmatrix}
[D_\mu^\nu] & 0 \\
[D_\mu^b] & [D_a^b]
    \end{pmatrix} \ ,
    \end{equation}
where $\{e_{I}^U\}$ is a \emph{specialized} frame in $T(U)$, $[D_\mu^\nu]$
and $[D_a^b]$ are non\ndash degenerate matrix\ndash valued functions on $U$,
and $D_\mu^a\colon U\to\field$.

    \begin{Defn}    \label{Defn6.1}
    The specialized frame $\{X_{I}\}$ over $U$ in $T(U)$, obtained
from~\eref{6.12} via an admissible transformation~\eref{3.4a} with matrix
\(
A = \Bigl(
    \begin{smallmatrix}
[D_\nu^\mu]^{-1}    & 0 \\
0           & [D_b^a]^{-1}
    \end{smallmatrix}
\Bigr) ,
\)
is called \emph{adapted to the frame $\{e_{I}\}$ for} $\Delta^h$.~%
\footnote{~%
Recall, here and below the adapted frames are defined only with respect to
frames $\{e_I\}=\{e_\mu,e_a\}$ such that $\{e_a\}$ is a basis for the
vertical distribution $\Delta^v$ over $U$, \ie $\{e_a|_p\}$ is a basis for
$\Delta_p^v$ for all $p\in U$. Since $\Delta^v$ is integrable, the relation
$e_a\in\Delta^v$ for all $a=n+1,\dots,n+r$ implies
$[e_a,e_b]_{\_}\in\Delta^v$ for all $a,b=n+1,\dots,n+r$.%
}
    \end{Defn}

	The frame $\{X_I\}$ adapted to $\{e_I\}$ is independent of the choice
of the specialized frame $\{e_I^U\}$ in~\eref{6.12} and can alternatively be
defined by
 $X_\mu=(\pi_*|_{\Delta^h})^{-1}\circ\pi_* (e_\mu)$ and $X_a=e_a$.

    If $\{u^I\}$ are bundle coordinates on $U$, the frame $\{X_I\}$
adapted to the coordinate frame $\bigl\{\frac{\pd}{\pd u^I}\bigr\}$ is said
to be \emph{adapted to the coordinates} $\{u^I\}$.

    According to~\eref{3.4a}, the adapted frame
$\{X_{I}\}=\{X_\mu,X_a\}$  is given by the equation
    \begin{align}   \label{6.13a}
(X_\mu,X_a)
& =
(e_\nu,e_b) \cdot
    \begin{bmatrix}
\delta_\mu^\nu  & 0 \\
+ \Gamma_\mu^b   &\delta_a^b
    \end{bmatrix}
=(e_\mu + \Gamma_\mu^b e_b , e_a)
    \end{align}
where the functions
$\Gamma_\mu^a\colon U\to\field$, called
\emph{(2\ndash index) coefficients of} $\Delta^h$ in $\{X_{I}\}$, are
defined by
    \begin{equation}    \label{6.14}
[\Gamma_\mu^a] := + [D_\nu^a]\cdot [D_\mu^\nu]^{-1} .
    \end{equation}

	A change $\{e_{I}\}\mapsto\{\tilde{e}_{I}\}$ with
    \begin{equation}    \label{6.15}
(\tilde{e}_\mu,\tilde{e}_a)
= (e_\nu,e_b) \cdot
    \begin{pmatrix}
[A_\mu^\nu] & 0 \\
[A_\mu^b]   & [A_a^b]
    \end{pmatrix}
=
(A_\mu^\nu e_\nu + A_\mu^b e_b , A_a^b e_b) ,
    \end{equation}
where $[A_\mu^\nu]$ and $[A_a^b]$ are non-degenerate matrix-valued functions
on $U$, which are constant on the fibres of $(E,\pi,M)$, and
$A_\mu^b\colon U\to\field$, entails the transformations
(see~\eref{6.12}--\eref{6.14})
    \begin{align}   \label{6.16}
(X_\mu,X_a) & \mapsto (\tilde{X}_\mu , \tilde{X}_a)
=
( \tilde{e}_\mu + \tilde{\Gamma}_\mu^b\tilde{e}_b,\tilde{e}_a )
=
( A_\mu^\nu X_\nu, A_a^b X_b)
=
(X_\nu,X_b) \cdot
    \begin{bmatrix}
A_\mu^\nu   & 0 \\
0   & A_a^b
    \end{bmatrix}
\\          \label{6.17}
\Gamma_\mu^a & \mapsto \tilde{\Gamma}_\mu^a
= \bigl([A_d^c]^{-1}\bigr)_b^a ( \Gamma_\nu^b A_\mu^\nu - A_\mu^b )
    \end{align}
of the frame $\{X_{I}\}$ adapted to $\{e_{I}\}$ and of the coefficients
$\Gamma_\mu^a$ of $\Delta^h$ in $\{X_{I}\}$, \ie  $\{\tilde{X}_{I}\}$
is the frame adapted to $\{\tilde{e}_{I}\}$ and $\tilde{\Gamma}_\mu^a$
are the coefficients of $\Delta^h$ in $\{\tilde{X}_{I}\}$.
	\begin{Note}	\label{Note6.1}
If $\{e_{I}\}$ and $\{\tilde{e}_{I}\}$ are adapted, then $A_\mu^b=0$.
	If $\{Y_I\}$ is a specialized frame, it is adapted to any frame
 $\{e_\mu=A_\mu^\nu Y_\nu,e_a=A_a^b Y_b\}$ and hence any specialized frame
can be considered as an adapted one; in particular, any specialized frame is
a frame adapted to itself. Obviously, the coefficients of a connection
identically vanish in a given specialized frame considered as an adapted one.
This leads to the concept of a \emph{normal frame}
to which is devoted the present paper.
Besides, from the above observation follows that the set of  adapted frames
coincides with the one of specialized frames.
	\end{Note}

    In particular, if $\{u^{I}\}$ and $\{\Tilde{u}^{I}\}$ are local
bundle coordinates with non\ndash empty intersection of their domains, we can
set
    \begin{equation}    \label{6.170}
e_I=\frac{\pd}{\pd u^I}
\quad
\Tilde{e}_I=\frac{\pd}{\pd \tilde{u}^I},
    \end{equation}
which entails
    \begin{equation}    \label{6.171}
A_{\mu}^{\nu} = \frac{\pd u^\nu}{\pd \tilde{u}^\mu}
\quad
A_{\mu}^{b} = \frac{\pd u^b}{\pd \tilde{u}^\mu}
\quad
A_{a}^{b} = \frac{\pd u^b}{\pd \tilde{u}^a}.
    \end{equation}
So, when the holonomic choice~\eref{6.170} is made, the
transformation~\eref{6.17} reduces to
    \begin{equation}    \label{3.22}
\Gamma_\mu^a \mapsto \tilde{\Gamma}_\mu^a
=
\Bigl( \frac{\pd \tilde{u}^a}{\pd u^b} \Gamma_\nu^b +
       \frac{\pd \tilde{u}^a}{\pd u^\nu}
\Bigr)  \frac{\pd u^\nu}{\pd \tilde{u}^\mu} .
    \end{equation}

    Let $(E,\pi,M)$ be a \emph{vector} bundle. According to the above-said in
this section, any \emph{adapted} frame $\{X_{I}\}=\{X_\mu,X_a\}$ in $T(E)$
is equivalent to a pair of frames in $T(M)$ and $E$ according to
    \begin{equation}    \label{6.22}
\{X_\mu,X_a\} \leftrightarrow
( \{E_\mu=\pi_*|_{\Delta^h}(X_\mu)\} , \{E_a=v^{-1}(X_a)\} ) .
    \end{equation}

    Suppose $\{X_{I}\}$ and $\{\tilde{X}_{I}\}$ are two adapted frames.
Then they are connected by (cf.~\eref{6.3a} and~\eref{6.16})
    \begin{equation}    \label{6.26}
\tilde{X}_\mu = (B_\mu^\nu\circ\pi) X_\nu
\quad
\tilde{X}_a = (B_a^b\circ\pi) X_b ,
    \end{equation}
where $[B_\mu^\nu]$ and $[B_a^b]$ as some non-degenerate matrix-valued
functions on $M$. The pairs of frames corresponding to them, in
accordance with~\eref{6.22}, are related via
    \begin{equation}    \label{6.27}
\tilde{E}_\mu = B_\mu^\nu E_\nu
\quad
\tilde{E}^a = B_a^b E_b
    \end{equation}
and \emph{vice versa}.

    \begin{Prop}    \label{Prop6.1}
Let $\Delta^h$ be a \emph{linear} connection on a \emph{vector} bundle
$(E,\pi,M)$~%
\footnote{~%
A connection on a vector bundle is linear if the generated by it parallel
transport is a linear mapping~\cite{Rahula}.%
}
and $\{X_\mu\}$ be the frame adapted for $\Delta^h$ to a frame $\{e_{I}\}$
such that $\{e_a\}$ is a basis for $\Delta^v$ and
    \begin{equation}    \label{6.28}
(e_\mu,e_a)|_U
=
(\pd_\nu,\pd_b)
\cdot
    \begin{bmatrix}
B_\mu^\nu\circ\pi   & 0 \\
(B_{c\mu}^{b}\circ\pi)\cdot E^c & B_a^b\circ \pi
    \end{bmatrix}
=
\bigl(
(B_\mu^\nu\circ\pi)\pd_\nu + ((B_{c\mu}^{b}\circ\pi)\cdot E^c)\pd_b ,
 (B_a^b\circ\pi)\pd_b
\bigl) ,
    \end{equation}
where $\pd_{I}:=\frac{\pd}{\pd u^{I}}$ for some local bundle
coordinates $\{u^{I}\}=\{u^\mu=x^\mu\circ\pi,u^b=E^b\}$ on $U\subseteq E$,
$[B_\mu^\nu]$ and $[B_a^b]$ are non\ndash degenerate matrix\ndash valued
functions on $U$, $B_{c\mu}^{b}\colon U\to\field$, and $\{E^a\}$ is the
coframe dual to $\{E_a=v^{-1}(X_a)\}$. Then the 2\ndash index coefficients
$\Gamma_\mu^a$ of $\Delta^h$ in $\{X_{I}\}$ have the representation
    \begin{equation}    \label{6.29}
\Gamma_\mu^a = - (\Gamma_{b\mu}^{a}\circ\pi)\cdot E^b
    \end{equation}
on $U$ for some functions $\Gamma_{b\mu}^{a}\colon U\to\field$, called
3\ndash index coefficients of $\Delta^h$ in $\{X_{I}\}$.
    \end{Prop}

    \begin{Rem} \label{Rem6.1}
The representation~\eref{6.29} is not valid for frames more general than the
ones given by~\eref{6.28}. Precisely, equation~\eref{6.29} is valid if and
only if ~\eref{6.28} holds for some local coordinates $\{u^{I}\}$ on $U$
--- see~\eref{6.17}.
    \end{Rem}

    \begin{Proof}
If $e_I=\frac{\pd}{\pd u^I}$ for some local coordinates $\{u^I\}$ on $E$, the
proposition coincides with the theorem in~\cite[p.~27]{Rahula}.
Writing~\eref{6.17} for the transformation
$\{\pd_{I}\}\mapsto\{e_{I}\}$, with $\{e_{I}\}$ given
by~\eref{6.28}, we get~\eref{6.29} with
\[
\Gamma_{b\mu}^{a}
= ( [B_d^e]^{-1})_c^a
  ( \lindex[\Gamma]{}{\pd}{}_{b\nu}^{c} B_\mu^\nu + B_{b\mu}^{c} ) ,
\]
where $\lindex[\Gamma]{}{\pd}{}_{b\nu}^{c}$ are the 3-index coefficients of
$\Delta^h$ in the frame adapted to the coordinates $\{u^{I}\}$.
    \end{Proof}

    Let $\{X_{I}\}$ and $\{\tilde{X}_{I}\}$ be frames adapted to
$\{e_{I}\}$ and $\{\tilde{e}_{I}\}$, respectively, with
(cf.~\eref{6.28})
    \begin{equation}    \label{6.30}
(\tilde{e}_\mu,\tilde{e}_a)
=
(e_\nu,e_b) \cdot
    \begin{bmatrix}
B_\mu^\nu\circ\pi   & 0 \\
(B_{c\mu}^{b}\circ\pi)\cdot E^c & B_a^b\circ \pi
    \end{bmatrix} \ ,
    \end{equation}
in which $\Delta^h$ admits 3-index coefficients. Then, due to~\eref{6.17}
and~\eref{6.29}, the 3\ndash index coefficients $\Gamma_{b\mu}^{a}$ and
$\tilde{\Gamma}_{b\mu}^{a}$ of $\Delta^h$ in  respectively $\{X_{I}\}$ and
$\{\tilde{X}_{I}\}$ are connected by
    \begin{equation}    \label{6.31}
\tilde{\Gamma}_{b\mu}^{a}
= \bigl([B_f^e]^{-1}\bigr)_c^a
  (\Gamma_{d\nu}^{c} B_\mu^\nu + B_{d\mu}^{c}) B_b^d .
    \end{equation}
It can easily be checked that the transformation
$\{e_{I}\}\mapsto\{\tilde{e}_{I}\}$, with $\{\tilde{e}_{I}\}$ given
by~\eref{6.30}, is the most general one that preserves the existence of
3\ndash index coefficients of $\Delta^h$ provided they exist in
$\{e_{I}\}$. Introducing the matrices
$\Gamma_\mu:=[\Gamma_{b\mu}^{a}]_{a,b=n+1}^{n+r}$,
$\tilde{\Gamma}_\mu:=[\tilde{\Gamma}_{b\mu}^{a}]_{a,b=n+1}^{n+r}$,
$B:=[B_b^a]$, and $B_\mu:=[B_{b\mu}^{a}]$, we rewrite~\eref{6.31} as
    \begin{equation}
    \tag{\protect\ref{6.31}$^\prime$}   \label{6.31'}
\tilde{\Gamma}_\mu
=
B^{-1}\cdot (\Gamma_\nu B_\mu^\nu +B_\mu)\cdot B.
    \end{equation}
A little below (see the text after equation~\eref{6.33}), we shall prove that
the compatibility of the developed formalism with the theory of covariant
derivatives requires further restrictions on the general transformed
frames~\eref{6.15} to the ones given by~\eref{6.30} with
    \begin{equation}    \label{6.32}
B_\mu= \tilde{E}_\mu(B)\cdot B^{-1} = B_\mu^\nu E_\nu(B) \cdot B^{-1} ,
    \end{equation}
where
\(
\tilde{E}_\mu
:=\pi_*|_{\Delta^h}(\tilde{X}_\mu)
=\pi_*|_{\Delta^h}((B_\mu^\nu\circ\pi)X_\nu)
= B_\mu^\nu E_\nu .
\)
In this case,~\eref{6.31'} reduces to
    \begin{equation}    \label{6.33}
\tilde{\Gamma}_\mu
= B_\mu^\nu B^{-1}\cdot (\Gamma_\nu\cdot B + E_\nu(B))
= B_\mu^\nu(B^{-1}\cdot\Gamma_\nu - E_\nu(B^{-1})) \cdot B .
    \end{equation}

    At last, a few words on the covariant derivatives operators $\nabla$
are in order. Without lost of generality, we define such an
operator
    \begin{equation}    \label{4.35}
    \begin{split}
\nabla & \colon\mathcal{X}(M)\times\Sec^1(E,\pi,M) \to  \Sec^0(E,\pi,M)
\\
\nabla & \colon (F,Y)\mapsto\nabla_FY
    \end{split}
    \end{equation}
via the equations
    \begin{subequations}    \label{4.40}
    \begin{align}   \label{4.40a}
\nabla_{F+G}Y &= \nabla_{F}Y + \nabla_{G}Y
\\          \label{4.40b}
\nabla_{fF}Y &= f\nabla_{F}Y
\\          \label{4.40c}
\nabla_{F}(Y+Z) &= \nabla_{F}Y + \nabla_{F}Z
\\          \label{4.40d}
\nabla_{F}(fY) &= F(f)\cdot Y + f\cdot \nabla_{F}Y ,
    \end{align}
    \end{subequations}
where
$F,G\in\mathcal{X}(M)$, $Y,Z\in\Sec^1(E,\pi,M)$, and $f\colon M\to\field$ is
a $C^1$ function.
Suppose $\{E_\mu\}$ is a
basis for $\mathcal{X}(M)$ and $\{E_a\}$ is a one for $\Sec^1(E,\pi,M)$.
Define the \emph{components} $\Gamma_{b\mu}^{a}\colon M\to\field$ of $\nabla$
in the pair of frames $(\{E_\mu\},\{E_a\})$ by
    \begin{equation}    \label{6.35}
\nabla_{E_\mu}(E_b) = \Gamma_{b\mu}^{a} E_a .
    \end{equation}
Then~\eref{4.40} imply
\[
\nabla_F Y = F^\mu( E_\mu(Y^a) +\Gamma_{b\mu}^{a} Y^b) E_a
\]
for $F=F^\mu E_\mu\in\mathcal{X}(M)$ and $Y=Y^aE_a\in\Sec^1(E,\pi,M)$. A
change $(\{E_\mu\},\{E_a\})\mapsto(\{\tilde{E}_\mu\},\{\tilde{E}_a\})$, given
via~\eref{6.27}, entails
    \begin{equation}    \label{6.36}
\Gamma_{b\mu}^{a}\mapsto \tilde{\Gamma}_{b\mu}^{a}
=
B_\mu^\nu \bigl([B_f^e]^{-1}\bigr)_c^a
                (\Gamma_{d\nu}^{c} B_b^d + E_\nu(B_b^c)) ,
    \end{equation}
as a result of~\eref{6.35}. In a more compact matrix form, the last result
reads
    \begin{equation}
    \tag{\protect\ref{6.36}$^\prime$}   \label{6.36'}
\tilde{\Gamma}_\mu = B_\mu^\nu B^{-1}\cdot (\Gamma_\nu\cdot B + E_\nu(B))
    \end{equation}
with $\Gamma_\mu:=[\Gamma_{_b\mu}^{a}]$,
$\tilde{\Gamma}_\mu:=[\tilde{\Gamma}_{b\mu}^{a}]$, and $B:=[B_b^a]$.

    Thus, if we identify the 3-index coefficients of $\Delta^h$, defined
by~\eref{6.29}, with the components of $\nabla$, defined by~\eref{6.35},~%
\footnote{~%
Such an identification is justified by the definition of $\nabla$ via the
parallel transport assigned to $\Delta^h$ or via a projection, generated by
$\Delta^h$, of a suitable Lie derivative on $\mathfrak{X}(E)$
--- see~\cite{Rahula
}.%
}
then the quantities~\eref{6.31'} and~\eref{6.36'} must coincide, which
immediately leads to the equality~\eref{6.32}. Therefore
    \begin{equation}    \label{6.37}
(e_\mu,e_a) \mapsto (\tilde{e}_\mu,\tilde{e}_a)
=
(e_\nu,e_b) \cdot
    \begin{bmatrix}
B_\mu^\nu\circ\pi & 0 \\
\bigl( (B_\mu^\nu E_\nu(B_d^b) (B^{-1})_c^d)\circ\pi \bigr) E^c    &
B_a^b\circ\pi
    \end{bmatrix}
\bigg|_{B=[B_a^b]}
    \end{equation}
is the most general transformation between frames in $T(E)$ such that the
frames adapted to them are compatible with the linear connection and
the covariant derivative corresponding to it. In particular, such are all
frames
 $\bigl\{\frac{\pd}{\pd u^{I}}\bigr\}$ in $T(E)$ induced by some vector
bundle coordinates $\{u^{I}\}$ on $E$ as the vector fibre coordinates
transform in a linear way like
\(
u^a \mapsto \tilde{u}^a=(B_b^a\circ\pi)\cdot u^b ;
\)
the rest members of the class of frames mentioned are obtained from them
via~\eref{6.37} with $e_{I}=\frac{\pd}{\pd u^{I}}$ and some non\ndash
degenerate matrix\ndash valued functions $[B_\mu^\nu]$ and $B$.

    If $\{X_{I}\}$ (resp.\ $\{\tilde{X}_{I}\}$) is the frame
adapted to a frame $\{e_{I}\}$ (resp.\ $\{\tilde{e}_{I}\}$), then
the change $\{e_{I}\} \mapsto\{\tilde{e}_{I}\}$, given
by~\eref{6.37}, entails $\{X_{I}\} \mapsto\{\tilde{X}_{I}\}$
with $\{\tilde{X}_{I}\}$ given by~\eref{6.26} (see~\eref{6.15}
and~\eref{6.16}). Since the last transformation is tantamount to
the change
    \begin{equation}    \label{6.38}
( \{E_\mu\}, \{E_a\} ) \mapsto ( \{\tilde{E}_\mu\},\{\tilde{E}_a\} )
    \end{equation}
of the basis of $\mathcal{X}(M)\times\Sec(E,\pi,M)$ corresponding
to $\{X_{I}\}$ via ~\eref{6.22}, \eref{6.26}, and~\eref{6.27}),
we can say that
the transition~\eref{6.38} induces the change~\eref{6.36} of the
3\ndash index coefficients of the connection $\Delta^h$. Exactly
the same is the situation one meets in the
literature~\cite{K&N-1,Warner,Poor} when covariant derivatives are
considered (and identified with connections).

    Regardless that the change~\eref{6.37} of the frames in $T(E)$ looks
quite special, it is the most general one that, through~\eref{6.16}
and~\eref{6.22}, is equivalent to an arbitrary change~\eref{6.38} of a basis
in $\mathcal{X}(M)\times\Sec(E,\pi,M)$, \ie of a pair of frames in $T(M)$ and
$E$.


\section {Normal frames: general case}
\label{Sect7.1}

    In the theory of linear connections on a manifold, the normal frames
are defined as frames in the tangent bundle space in which the connections'
(3\ndash index) coefficients vanish on some subset of the
manifold~\cite{K&N-1,Poor,ORai,bp-Frames-n+point,bp-Frames-path,
bp-Frames-general}. The definition of normal frames for a connection on a
vector bundle is practically the same, the only difference being that these
frames are in the bundle space, not in the tangent bundle space over the base
space~\cite{bp-NF-D+EP}. The present section is devoted to the introduction
of normal frames for general connections on fibre bundles and some their
properties.

	To save some space and for brevity, in what follows we shall not
indicate explicitly that the frames $\{e_I\}=\{e_\mu,e_a\}$, with respect to
which the adapted frames are defined, are such that $\{e_a\}$ is a (local)
basis for the vertical distribution $\Delta^v$ on the bundle considered.

    \begin{Defn}    \label{Defn7.1}
Given a connection $\Delta^h$ on a bundle $(E,\pi,M)$ and a subset
$U\subseteq E$. A
\emph{frame $\{X_{I}\}$ in $T(E)$ adapted to a frame $\{e_{I}\}$ in}
 $T(E)$ and defined over an open subset $V$ of $E$ containing or equal to
$U$, $V\supseteq U$, is called
\emph{normal for $\Delta^h$ over/on $U$ (relative to $\{e_I\}$)}
if in it all (2\ndash index) coefficients  $\Gamma_\mu^a$ of $\Delta^h$
vanish everywhere on $U$. Respectively, $\{X_{I}\}$ is \emph{normal for
$\Delta^h$ along a mapping} $g\colon Q\to E$, $Q\not=\varnothing$, if
$\{X_{I}\}$ is normal for $\Delta^h$ over the set $g(Q)$.
    \end{Defn}

    Let $\{X_{I}\}$ be the frame in $T(E)$ adapted to a frame
$\{e_{I}\}$ in $T(E)$ over an open subset $V\subseteq E$. Then the
frame $\{\tilde{X}_{I}\}$ in $T(E)$ adapted to a frame
$\{\tilde{e}_{I}\}$, given by~\eref{6.15}, in $T(E)$ over the
same subset $V$ is normal for $\Delta^h$ over $U\subseteq V$ if and only if
    \begin{equation}    \label{7.1}
(A_\mu^\nu \Gamma_\nu^b - A_\mu^b)|_U = 0 ,
    \end{equation}
due to~\eref{6.16} and~\eref{6.17}. Since $\Gamma_\mu^b$ depend
only on $\Delta^h$ and $\{e_{I}\}$, the existence of solutions
of~\eref{7.1}, relative to $A_\mu^\nu$ and $A_\mu^b$, and their
properties are completely responsible for the existence and the
properties of frames normal for $\Delta^h$ over $U$. For that
reason, we call~\eref{7.1} the \emph{(system of) equation(s) of
the normal frames for $\Delta^h$ over $U$} or simply the
\emph{normal frame (system of) equation(s)} (for $\Delta^h$ over
$U$).

    In the most general case, when no additional restrictions on the
frames considered are imposed, the normal frames
equation~\eref{7.1} is a system of  $nr$ \emph{linear algebraic
equations} for $nr+n^2$ variables and, consequently, it has a
solution depending on $n^2$ independent parameters. In particular,
if we choose the functions $A_\mu^\nu\colon U\to\field$ (with
$\det[A_\mu^\nu]\not=0,\infty$) as such parameters, we can write
the general solution of~\eref{7.1} as
    \begin{equation}    \label{7.2}
( \{A_\mu^\nu\} , \{A_\mu^b\} ) |_U
=
( \{A_\mu^\nu\} , \{\Gamma_\nu^b A_\mu^\nu\} )|_U .
    \end{equation}

	It should be noted, equation~\eref{7.1} or its general
solution~\eref{7.2} defines the frame $\{\tilde{e}_I\}$ and the frame
$\{\tilde{X}_I\}$  adapted to $\{\tilde{e}_I\}$ only on $U$ and leaves them
completely arbitrary on $V\setminus U$, if it is not empty.

    \begin{Prop}    \label{Prop7.1}
Let $\{X_{I}\}$  be the frame adapted to a frame $\{e_I\}$ in
$T(V)\subseteq T(E)$ defined over an open set $V\subseteq E$ and
$\Gamma_{\mu}^{a}$ be the coefficients of a connection $\Delta^h$ in
$\{X_{I}\}$. Then all frames $\{\tilde{X}_{I}\}$ normal on $U\subseteq V$
for the connection $\Delta^h$ are adapted to frames $\{\tilde{e}_{I}\}$ given
on $U$ by
    \begin{equation}    \label{7.4}
\tilde{e}_\mu |_U = (A_\mu^\nu(e_\nu + \Gamma_\nu^b e_b))|_U
\quad
\tilde{e}_a|_U = (A_a^b e_b)|_U .
    \end{equation}
where $[A_\mu^\nu]$ and $[A_a^b]$ are non-degenerate matrix-valued
functions on $V$ which are constant on the fibres of $(E,\pi,M)$. Moreover,
the frame $\{\tilde{X}_I\}$ adapted on $V$ to $\{\tilde{e}_I\}$, given
by~\eref{7.4} (and hence normal on $U$), is such that
    \begin{equation}    \label{7.3}
\tilde{X}_\mu|_U = (A_\mu^\nu X_\nu)|_U = \tilde{e}_\mu |_U
\quad
\tilde{X}_a|_U = (A_a^b X_b)|_U = \tilde{e}_a |_U  ,
    \end{equation}
    \end{Prop}

    \begin{Proof}
Apply~\eref{6.16},~\eref{6.15}, and~\eref{6.13a} for the
choice~\eref{7.2}.
    \end{Proof}

    The equations~\eref{7.3} are not accidental as it is stated by the
following assertion.

    \begin{Prop}    \label{Prop7.2}
The frame $\{\tilde{X}_{I}\}$ in $T(E)$ adapted to a frame
$\{\tilde{e}_{I}\}$ in $T(E)$ and defined over an open set $V\subseteq E$ is
normal on $U\subseteq V$ if and only if on $U$ is fulfilled
    \begin{equation}    \label{7.5}
\tilde{X}_{I}|_U = \tilde{e}_{I}|_U .
    \end{equation}
    \end{Prop}

    \begin{Proof}
Apply~\eref{6.13a} or~\eref{6.16} and definition~\ref{Defn7.1}.
    \end{Proof}

    Thus one can equivalently \emph{define the normal frames as adapted
frames that coincide on some set with the frames they are adapted to} or as
frames (in the tangent bundles space over the bundle space) that coincide
on some set with the frames adapted to them.

	Since any specialized frame is adapted to itself (see
definition~\ref{Defn6.1} and~\eref{6.12}, with $D_I^J=\delta_I^J$),
\emph{the sets of normal, specialized, and adapted frames are identical}.

    As we see from proposition~\ref{Prop7.1}, which gives a complete
description of the normal frames, the theory of normal frames in
the most general setting is trivial. It becomes more interesting
and richer if the class of frames $\{e_{I}\}$, with respect to
which are defined the adapted frames, is restricted in one or
other way. To the theory of normal frames, adapted to such
restricted classes of frames in $T(E)$, are devoted the next two sections.

\section{Normal frames adapted to holonomic frames}
    \label{Sect7.2}

    The most natural class of frames in $T(E)$ relative to which the
adapted, in particular normal, frames are defined is the one of
holonomic frames induced by local coordinates on $E$ (see
Subsect.~\ref{Subsect3.2}). To specify the consideration of the
previous section to normal frames adapted to local coordinates
on $E$, we set $e_{I}=\frac{\pd}{\pd u^{I}}$ and
$\tilde{e}_{I}=\frac{\pd}{\pd \tilde{u}^{I}}$, where
$\{u^{I}\}$ and $\{\tilde{u}^{I}\}$ are local coordinates on
$E$ whose domains have a non\ndash empty intersection $V$ and
$U\subseteq V$. Then the matrix $[A_{I}^{J}]$ in~\eref{7.1} is
given by~\eref{3.3} (as $\{e_{I}\}\mapsto\{\tilde{e}_{I}\}$
reduces to~\eref{3.2a}), so that the normal frame
equation~\eref{7.1} reduces to the \emph{normal coordinates equation} (see
also~\eref{3.22})
    \begin{equation}    \label{7.6}
\Bigl(
\frac{\pd\tilde{u}^a}{\pd u^b} \Gamma_\mu^b
+ \frac{\pd\tilde{u}^a}{\pd u^\mu} \Bigr)
\Big|_U
=0 ,
    \end{equation}
due to~\eref{3.1}, which is a first order system of $nr$ \emph{linear
partial differential equations} on $U$ relative to the $r$ unknown functions
$\{\tilde{u}^{n+1},\dots,\tilde{u}^{n+r}\}$.

    Since the connection $\Delta^h$ is supposed given and fixed, such are its
coefficients $\Gamma_\mu^b$ in $\bigl\{\frac{\pd}{\pd u^I}\bigr\}$.
Therefore the existence, uniqueness and other properties of the solutions
of~\eref{7.6} strongly depend on the set $U$ (which is in the intersection
of the domains of the local coordinates $\{u^{I}\}$ and $\{\tilde{u}^{I}\}$
on $E$).

    \begin{Prop}    \label{Prop7.3}
If the normal frame equation~\eref{7.6} has solutions, then all frames
$\{\Tilde{X}_I\}$ normal on $U\subseteq E$  and adapted to local
coordinates, defined on an open set $V\subseteq E$ such that $V\supseteq U$,
are described by
    \begin{equation}    \label{7.7}
\tilde{X}_\mu |_U
= (A_\mu^\nu X_\nu)|_U = \frac{\pd}{\pd\tilde{u}^\mu} \Big|_U
\quad
\tilde{X}_a|_U
= (A_a^b X_b)|_U = \frac{\pd}{\pd\tilde{u}^a} \Big|_U ,
    \end{equation}
where $\{X_{I}\}$ is the frame adapted to some arbitrarily fixed local
coordinates $\{u^{I}\}$, defined on an open set containing or equal to $V$,
$\{\Tilde{u}^I\}$ are local coordinates with domain $V$ and such that
$\Tilde{u}^a$ are solutions of~\eref{7.6}, and $A_{I}^{J}=\frac{\pd
u^{J}}{\pd\tilde{u}^{I}}$ on the intersection of the domains of $\{u^{I}\}$
and $\{\tilde{u}^{I}\}$.
    \end{Prop}

    \begin{Proof}
Apply proposition~\ref{Prop7.1} for $e_{I}=\frac{\pd}{\pd u^{I}}$ and
$\tilde{e}_{I}=\frac{\pd}{\pd \tilde{u}^{I}}$ and then use~\eref{3.2a}
and~\eref{3.3}.
    \end{Proof}

    This simple result gives a complete description of all normal frames, if
any, adapted to (local) holonomic frames. It should be understood clearly,
normal on $U$ is the frame $\{\Tilde{X}_I\}$, adapted to
$\bigl\{\frac{\pd}{\pd\Tilde{u}^I}\bigr\}$ and coinciding with it on $U$,
but not the frame $\bigl\{\frac{\pd}{\pd\Tilde{u}^I}\bigr\}$; in particular,
the frame $\bigl\{\frac{\pd}{\pd\Tilde{u}^I}\bigr\}$ is holonomic while the
frame $\{\Tilde{X}_I\}$ need not to be holonomic, even on $U$, if the
connection considered does not satisfies some additional conditions, like the
vanishment of its curvature on $U$.

    Consider now briefly the existence problem for the solutions
of~\eref{7.6}. To begin with, we emphasize that in~\eref{7.6} enter only the
fibre coordinates $\{\tilde{u}^a\}$, so that it leaves the basic ones
$\{\tilde{u}^\mu\}$ completely arbitrary.

    \begin{Prop}    \label{Prop7.4}
If $E$ is of class $C^2$, $p\in E$ is fixed, and $U=\{p\}$, then the general
solution of~\eref{7.6} is
    \begin{equation}    \label{7.8}
\tilde{u}^a(q)
=
g^a
+ g_b^a\{ - \Gamma_\mu^b(p) (q^\mu-p^\mu) + (q^b-p^b) \}
+ f_{{I}{J}}^{a}(q) (q^{I}-p^{I})(q^{J}-p^{J}),
    \end{equation}
where $g^a$ and $g_b^a$ are constants in $\field=\field[R],\field[C]$,
$\det[g_b^a]\not=0,\infty$, the point  $q$ is in the domain $V$ of
$\{u^{I}\}$, $q^{I}:=u^{I}(q)$, $p^{I}:=u^{I}(p)$, and
$f_{{I}{J}}^{a}$ are $C^2$ functions on $V$ such that they and their
first partial derivatives are bounded when $q^I\to p^I$.
    \end{Prop}

    \begin{Proof}
Expand
\( \tilde{u}^a(q)
=f^a(u^1(q),\dots,u^n(q),\dots,u^{n+r}(q))=f^a(q^1,\dots,q^{n+r})
\)
into a Taylor's first order polynomial with remainder term
quadratic in $(q^{I}-p^{I})$ and insert the result
into~\eref{7.6}. In this way, we get~\eref{7.8} with
$g^a=\tilde{u}^a(p)$ and
$g_b^a=\frac{\pd\tilde{u}^a}{\pd u^b}\big|_{p}$.
    \end{Proof}

	Now we would like to investigate the existence of solutions
of~\eref{7.6} along paths $\beta \colon J\to E$, \ie for $U=\beta(J)$. The main
result is formulated below as proposition~\ref{Prop7.5}. For its proof, we
shall need the following lemma.

    \begin{Lem}    \label{2-Lem6.1-0}
Let $\gamma \colon J\to M$ be a regular $C^1$ injective path in a $C^3$ real
manifold $M$. For every $s_0\in J$, there exists a chart $(U_1,x)$ of $M$ such
that $\gamma(s_0)\in U_1$ and $x(\gamma(s))=(s,\bs{t}_0)$ for all $s\in J$ such
that $\gamma(s)\in U_1$ and some fixed $\bs{t}_0\in\field[R]^{\dimR M-1}$.
    \end{Lem}

    \begin{Proof}
	Let $s_0\in J$ be a point in $J$ which is not an end point of $J$, if
any, and $(U,y)$ be a chart with $\gamma(s_0)$ in its domain,
$U\ni\gamma(s_0)$, and $y\colon U\to\field[R]^{\dimR M}$. From the regularity
of $\gamma$, $\dot\gamma\not=0$, follows that at least one of the numbers
$\dot\gamma_y^1(s_0),\ldots,\dot\gamma_y^{\dimR M}(s_0)$, where
$\gamma_{y}^{i}:=y^i\circ\gamma$, is non\ndash zero.
We, without lost of generality, choose this non\ndash vanishing component to
be $\dot\gamma_y^{1}(s_0)$.%
\footnote{\label{2-renumbering}%
If it happens that $\dot\gamma_y^{1}(s_0)=0$ and
$\dot\gamma_y^{i_0}(s_0)\not=0$ for some $i_0\not=1$, we have simply to
renumber the local coordinates to get $\dot\gamma_y^{1}(s_0)\not=0$.
Practically this is a transition to new coordinates $\{y^i\}\to\{z^i\}$ with
$z^1=y^{i_0}$ and, for instance, $z^{i_0}=y^{1}$ and $z^{i}=y^{i}$ for
$i\not=1,i_0$, in which the first component of $\dot\gamma$ is non-zero. We
suppose that, if required, this coordinate change is already done. If
occasionally it happens that $\dot\gamma_y^{j_0}(s)\not=0$ for all $s\in J$
and fixed $j_0$, it is extremely convenient to take this particular component
of $\dot\gamma$ as $\dot\gamma_y^1$ --- see the next sentence.%
}
Then, due to the continuity of $\dot\gamma$ ($\gamma$ is of class $C^1$) and
according to the implicit function
theorem~\cite[chapter~III, \S~8]{Schwartz/Analysis-1},
\cite[sect.~1.37 and~1.38]{Warner},
\cite[chapter~10, sect.~2]{Dieudonne},
there exists an \emph{open} subinterval $J_1\subseteq J$ containing $s_0$,
$J_1\ni s_0$, and such that $\dot\gamma^1|_{J_1}\not=0$ and the restricted
mapping $\gamma_{y}^{1}|_{J_1}\colon J_1\to \gamma_{y}^{1}(J_1)$ is a $C^1$
diffeomorphism on its image. Define a neighborhood
\[
U_1 := \bigl\{ p | p\in U,\ y^1(p)\in\gamma_{y}^{1}(J_1) \bigr\}
    = y^{-1}\bigl( \gamma_y^1(J_1)\times\field[R]^{\dimR M -1} \bigr)
\ni\gamma(s_0)
\]
and a chart $(U_1,x)$ with local coordinate functions
	\begin{equation}	\label{2-6.20}
	\begin{split}
x^1 &:= \bigl(\gamma_{y}^{1}|_{J_1}\bigr)^{-1} \circ y^1
\\
x^k &:= y^k - \gamma_y^k\circ x^1 + t_{0}^{k}
\qquad
k=2,\ldots,\dimR M
	\end{split}
	\end{equation}
where $t_{0}^{k}\in\field[R]$ are constant numbers.
Since
$\frac{\pd x^1}{\pd y^j}=\frac{1}{\dot\gamma_y^1}\delta_{j}^{1}$,
$\frac{\pd x^k}{\pd y^1}=-\frac{\dot\gamma_y^k}{\dot\gamma_y^1}$,
for $k\ge2$, and
$\frac{\pd x^k}{\pd y^l}=\delta_{l}^{k}$ for $k,l\ge2$,
the Jacobian of the change $\{y^i\}\to\{x^i\}$ at $p\in U_1$ is
$\frac{1}{\dot\gamma^1(p)}\not=0,\infty$. Consequently
 $x\colon U_1\to J_1\times\field[R]^{\dimR M-1}$ is really a coordinate
homeomorphism with coordinate functions $x^i$.

	In the new chart $(U_1,x$),
the coordinates of $\gamma(s)$, $s\in J_1$ are
	\begin{equation}	\label{2-6.20gamma}
\gamma^1(s) := (x^1\circ\gamma)(s) = s,
\quad
\gamma^k(s) := (x^k\circ\gamma)(s) = t_{0}^{k},\ k\ge2 ,
	\end{equation}
\ie $x(\gamma(s))=(s,\bs{t}_0)$ for some
$\bs{t}_0=(t_{0}^{2},\ldots,t_{0}^{\dimR M})\in\field[R]^{\dimR M-1}$.
    \end{Proof}

	Lemma~\ref{2-Lem6.1-0} means that the chart $(U_1,x)$ is so luckily
chosen that the first coordinate in it of a point along $\gamma$ coincides with
the value of the corresponding path's parameter, the other coordinates being
constant numbers. Moreover, in $U_1$ the path $\gamma$ can be considered as a
representative of a family of paths
$\eta(\cdot,\bs{t})\colon J_1\to M$,
$\bs{t}\in\field[R]^{\dimR M-1}$,
defined by $\eta(s,\bs{t}):=x^{-1}(s,\bs{t})$ for
 $(s,\bs{t})\in J_1\times\field[R]^{\dimR M-1}$; indeed,
 $\gamma=\eta(\cdot,\bs{t}_0)$
or
 $\gamma(s)=\eta(s,\bs{t}_0)$,  $s\in J_1\subseteq J$.

    \begin{Prop}    \label{Prop7.5}
Let $\Delta^h$ be a $C^1$ connection on a real $C^3$ bundle
$(E,\pi,M)$, $n=\dim M\ge1$, $r=\dim\pi^{-1}(x)\ge1$ for $x\in M$. Let
$\beta\colon J\to E$ be an injective regular  $C^1$ path such that its tangent
vector $\dot{\beta}(s)$ at $s$ is \emph{not} a vertical vector for all $s\in
J$, $\dot{\beta}(s)\not\in\Delta_{\beta(s)}^v$; in particular, the path $\beta$
can be horizontal, \ie $\dot{\beta}(s)\in\Delta_{\beta(s)}^h$ for all $s\in J$,
but generally the vector $\dot{\beta}(s)$ can have also and a vertical
component for some or all $s\in J$.
	Then, for every $s_0\in J$, there exist a neighborhood $U_1$ of the
point $\beta(s_0)$ in $E$ and bundle coordinates $\{\tilde{u}^I\}$ on $U_1$
which are solutions of~\eref{7.6} for $U=U_1\cap\beta(J)=\beta(J_1)$, with
$J_1:=\{s\in J : \beta(s)\in U_1\}$, \ie along the restricted path
$\beta|_{J_1}$.
	All such bundle coordinates $\{\tilde{u}^{I}\}$ are given via the
equation~\eref{7.10} below.
    \end{Prop}

    \begin{Proof}
Consider the chart $(U_1,u)$ with $U_1\ni\beta(s_0)$ provided by
lemma~\ref{2-Lem6.1-0} for $E$ and $\beta$ instead of $M$ and $\gamma$,
respectively. For any $p\in U_1$, there is a unique $(s,\bs{t})\in
J_1\times\field[R]^{\dimR E-1}$ such that $p=u^{-1}(s,\bs{t})$, i.e., in the
coordinates $\{u^I\}$ associated to $u$, the coordinates of $p$ are $u^1(p)=s$
and $u^I(p)=t^I\in\field[R]$ for $I\ge2$. Besides, we have
$u(\beta(s))=(s,\bs{t_0})$ for all $s\in J_1$ and some fixed
$\bs{t_0}\in\field[R]^{\dimR E-1}$.

	Since $\dot{\beta}(s)$ is not a vertical vector for all $s\in J$, the
coordinates $\{u^I\}$ can be chosen to be \emph{bundle} coordinates. For the
purpose, in the proof of lemma~\ref{2-Lem6.1-0} one must choose $\{y^I\}$ as
bundle coordinates and to take for $\dot{\beta}_y^1(s_0)$ any non\ndash
vanishing component between $\dot{\beta}_y^1(s_0),\dots,\dot{\beta}_y^n(s_0)$,
\viz if $\dot{\beta}_y^1(s_0)\not=0$ the proof goes as it is written and, if
$\dot{\beta}_y^1(s_0)=0$, choose some $\mu_0$ such that
$\dot{\beta}_y^{\mu_0}(s_0)\not=0$ and make, e.g., the change
$\dot{\beta}_y^1(s_0)\leftrightarrow\dot{\beta}_y^{\mu_0}(s_0)$. This, together
with~\eref{2-6.20}, with $u^I$ for $x^k$, ensures that $\{y^I\}\mapsto\{u^I\}$
is an admissible change, so that $\{u^I\}$ are bundle coordinates if the
initial coordinates $\{y^I\}$ are such ones.

	Let $\{u^I\}$ be so constructed bundle coordinates and $\eta:=u^{-1}$,
so that $\beta(s)=\eta(s,\bs{t}_0)$. Expanding $\tilde{u}^a(\eta(s,\bs{t}))$
into a first order Taylor's polynomial at the point $\bs{t}_0\in K$, we find
the general solution of~\eref{7.6}, with $U=\beta(J_1)=U_1\cap\beta(J)$, in
the form
    \begin{multline}    \label{7.10}
\tilde{u}^a(\eta(s,\bs{t}))
=
B^a(s)
+ B_b^a(s)\{ - \Gamma_\mu^b(\beta(s)) [ u^\mu(\eta(s,\bs{t}))-u^\mu(\beta(s))]
+[ u^b(\eta(s,\bs{t})) - u^{b}(\beta(s)) ] \}
\\
+B_{{I}{J}}^{a}(s,\bs{t};\eta)
[u^{I}(\eta(s,\bs{t})) - u^{I}(\beta(s))]
[u^{J}(\eta(s,\bs{t})) - u^{J}(\beta(s))] ,
    \end{multline}
where $B^a,B_b^a\colon J_1\to\field=\field[R]$, $\det[B_b^a]\not=0,\infty$,
and the $C^1$ functions $B_{{I}{J}}^{a}$ and their first
partial derivatives are bounded when $\bs{t}\to\bs{t}_0$. (Notice,
the terms with $\mu=1$ and/or ${I}=1$ and/or ${J}=1$ do not
contribute in~\eref{7.10} as $u^1(\eta(s,\bs{t}))\equiv s$ and,
besides, the functions $B_{{I}{J}}^{a}$ can be taken
symmetric in ${I}$ and ${J}$, $B_{{I}{J}}^{a}=B_{{J}{I}}^{a}$.)
    \end{Proof}

    \begin{Rem}    \label{Rem7.0}
If there is $s_0\in J$ for which $\dot{\beta}(s_0)$ is a vertical vector,
$\dot{\beta}(s_0)\in\Delta_{\beta(s_0)}^v$, then proposition~\ref{Prop7.5}
remains true with the only correction that the coordinates $\{u^I\}$ will
\emph{not be bundle} coordinates. If this is the case, the constructed
coordinates $\{\tilde{u}^I\}$ will be solutions of~\eref{7.6}, but we cannot
assert that they are bundle coordinates which are (locally) normal along
$\beta$ in a neighborhood of the point $\beta(s_0)$.
    \end{Rem}

	Proposition~\ref{Prop7.5} can be generalized by requiring $\beta$ to be
locally injective instead of injective, \ie for each $s\in J$ to exist a
subinterval $J_s\subseteq J$ such that $J_s\ni s$ and the restricted path
$\beta|_{J_s}$ to be injective.
	Besides, if one needs a version of the above results for complex
bundles, they should be considered as real ones (with doubled dimension of the
manifolds) for which are applicable the above considerations.

    \begin{Cor} \label{Cor7.1}
At any arbitrarily fixed point in $E$ and/or along a given injective regular
$C^1$ path in $E$, whose tangent vector is not vertical, there exist (possibly
local, in the latter case) normal frames.
    \end{Cor}

    \begin{Proof}
See definition~\ref{Defn7.1}, propositions~\ref{Prop7.4}
and~\ref{Prop7.5}, and equation~\eref{7.6}. If the path is not
contained in a single coordinate neighborhood, one should cover its
image in the bundle space with such neighborhoods and, then, to
apply proposition~\ref{Prop7.5}; in the intersection of the
coordinate domains, the uniqueness (and, possibly, continuity or
differentiability) of the normal frames may be lost.
    \end{Proof}

    \begin{Defn}    \label{Defn7.2}
Local bundle coordinates $\{\tilde{u}^{I}\}$, defined on an open set
$V\subseteq E$, will be called \emph{normal} on $U\subseteq V$ for a
connection $\Delta^h$ if the frame $\{\Tilde{X}_I\}$ in $T(E)$ adapted to
$\bigl\{ \frac{\pd}{\pd\tilde{u}^{I}} \bigr\}$ over $V$ is normal for
$\Delta^h$ on $U$.
    \end{Defn}

    Corollary~\ref{Cor7.1} implies the existence of coordinates normal at
a given point or (locally) along a given injective path whose tangent vector
is not vertical; in particular, there exist coordinates normal along an
injective horizontal path. However, normal coordinates generally do not exist
on more general subsets of the bundle space $E$. A criterion for existence of
coordinates normal on sufficiently general subsets $U\subseteq E$, e.g on
`horizontal' submanifolds, is given by theorem~\ref{ThmA.1} in the Appendix. In
particular, we have the following
corollary from this theorem.

    \begin{Prop}    \label{Prop7.6}
If $\Delta^h$ is a $C^1$ connection, $U$ is an open set in $E$, and
normal frames for $\Delta^h$ on $U$ exist, then there are holonomic such
frames if $\Delta^h$ is flat on $U$. Said otherwise, the system
of equations~\eref{7.6} may admit solutions on an open set $U$ if
    \begin{equation}    \label{7.9}
R_{\mu\nu}^{a}|_U=0
    \end{equation}
where
    \begin{align}   \label{3.24a}
R_{\mu\nu}^{a}
& = \pd_\mu(\Gamma_\nu^a) - \pd_\nu(\Gamma_\mu^a)
 + \Gamma_\mu^b\pd_b(\Gamma_\nu^a)
 - \Gamma_\nu^b\pd_b(\Gamma_\mu^a)
 = X_\mu(\Gamma_\nu^a) -  X_\nu(\Gamma_\mu^a) .
    \end{align}
are the (fibre) components of the curvature of $\Delta^h$ in some frame
$\{X_I\}$ on $E$ adapted to a holonomic one
(see also~\cite{Rahula
}%
).
    \end{Prop}


	\begin{Rem}	\label{Rem7.1}
However, in the general case the flatness of a connection on an open set is
only a necessary, but not sufficient, condition for the existence of
coordinates normal on that set --- see theorem~\ref{ThmA.1} in the Appendix.
Exceptiones are the linear connections on vector bundles --- see
remark~\ref{RemA.3} in the Appendix.
	One can easily show that part of the integrability conditions
for~\eref{7.6} for an open set $U$ are
    \begin{equation}    \label{7.9-1}
0
=
  \frac{\pd^2\tilde{u}^a}{\pd u^\nu \pd u^\mu}
- \frac{\pd^2\tilde{u}^a}{\pd u^\mu \pd u^\nu}
\equiv
\frac{\pd \tilde{u}^a}{\pd u^b} R_{\mu\nu}^{b}
    \end{equation}
from where proposition~\ref{Prop7.6} immediately follows. However, the
flatness of the connection on $U$ generally does not imply the rest of the
integrability conditions, \viz
\(
  \frac{\pd^2\tilde{u}^a}{\pd u^b \pd u^\mu}
- \frac{\pd^2\tilde{u}^a}{\pd u^\mu \pd u^b}
= 0
\)
and
\(
  \frac{\pd^2\tilde{u}^a}{\pd u^b \pd u^c}
- \frac{\pd^2\tilde{u}^a}{\pd u^c \pd u^b}
= 0 .
\)
	\end{Rem}

    The combination of propositions~\ref{Prop7.6} and~\ref{Prop7.2}
implies the non\ndash existence of coordinates normal on an open set
for non\ndash flat (non\ndash integrable) connections.

\section{Normal frames on vector bundles}
    \label{Sect7.3}

    The normal frames for covariant derivative operators (linear
connections), other derivations, and  linear transports along paths are known
and studied objects in vector bundles~\cite{bp-NF-D+EP,bp-NF-LTP}. The
goal of the present section is to be made a link between them and the
general theory of Sect.~\ref{Sect7.1}.

    Consider a linear connection $\Delta^h$ on a vector bundle $(E,\pi,M)$,
\ie a connection the assigned to which parallel transport is a linear
mapping. Let the frame $\{e_{I}\}$ in $T(E)$ be given by~\eref{6.28} and
$\{X_{I}\}$ be the frame adapted to $\{e_{I}\}$ for $\Delta^h$. Then, by
proposition~\ref{Prop6.1}, the 2- and 3\ndash index coefficients of
$\Delta^h$ are connected via~\eref{6.29} in which $\{E^a\}$ is a frame over
$M$ dual to the frame $\{E_a=v^{-1}(X_a)\}$ in $E$, with $v$ defined
by~\eref{4.1-1}.

    \begin{Prop}    \label{Prop7.7}
A frame $\{X_{I}\}$ is normal on $U\subseteq E$ for a linear connection
$\Delta^h$ if and only if in it vanish the 3\ndash index coefficients of
$\Delta^h$ on $\pi(U)\subseteq M$,
    \begin{equation}    \label{7.11}
\Gamma_\mu^a|_U = 0 \iff \Gamma_{b\mu}^{a}|_{\pi(U)} = 0 .
    \end{equation}
    \end{Prop}

    \begin{Proof}
Since the linear forms $E^{n+1}|_p,\dots,E^{n+r}|_p$ are linearly
independent for all $p\in U$, the assertion follows from
equation~\eref{6.29}.
    \end{Proof}

    Combining proposition~\ref{Prop7.7} with~\eref{6.31}, we see that the
normal frame equation~\eref{7.1} in vector bundle is equivalent to
    \begin{equation}    \label{7.12}
(B_\mu^\nu \Gamma_{b\nu}^{a} + B_{b\mu}^{a})|_{\pi(U)} = 0
    \end{equation}
or to its matrix variant (see also~\eref{6.31'};
$\Gamma_\nu:=[\Gamma_{b\nu}^{a}]$, $B_\mu:=[B_{b\mu}^{a}]$)
    \begin{equation}
    \tag{\protect\ref{7.12}$^\prime$}   \label{7.12'}
(B_\mu^\nu \Gamma_\nu + B_\mu)|_{\pi(U)} = 0 .
    \end{equation}
Taking into account~\eref{7.12} and~\eref{6.30}, we can assert that the frame
$\{\tilde{X}_{I}\}$ adapted to the frame
    \begin{equation}    \label{7.13}
(\tilde{e}_\mu,\tilde{e}_a)
=
(e_\nu,e_b)
\cdot
    \begin{bmatrix}
B_\mu^\nu\circ\pi   & 0 \\
-( (B_\mu^\lambda \Gamma_{c\lambda}^{b}) \circ\pi ) \cdot E^c
            & B_a^b\circ\pi
    \end{bmatrix} \ ,
    \end{equation}
where $[B_\mu^\nu]$ and $[B_b^a]$ are non-degenerate matrix-valued functions,
is normal on $U$ for $\Delta^h$ and hence
$\tilde{X}_{I}=\tilde{e}_{I}$, by virtue of proposition~\ref{Prop7.2}.
Recall (see~\eref{6.15}, \eref{6.16}, and~\eref{6.27}), the change
$\{e_{I}\}\mapsto\{\tilde{e}_{I}\}$, given by~\eref{7.13}, entails
$\{X_{I}\}\mapsto\{\tilde{X}_{I}\}$, where
    \begin{equation}    \label{7.14}
\tilde{X}_\mu = (B_\mu^\nu\circ\pi) X_\nu
\quad \tilde{X}_a = (B_a^b\circ\pi) X_b,
    \end{equation}
which is equivalent to
$\{E_{I}\}\mapsto\{\tilde{E}_{I}\}$ with
    \begin{equation}    \label{7.15}
\tilde{E}_\mu = B_\mu^\nu E_\nu
\quad \tilde{E}_a = B_a^b E_b .
    \end{equation}
Here (see~\eref{6.22}) $\{E_\mu=\pi_*|_{\Delta^h}(X_\mu)\}$ is a frame in
$T(M)$ and $\{E_a=v^{-1}(X_a)\}$ is a frame in $E$.

    Thus, if additional restriction are not imposed, the theory of normal
frames in vector bundles is rather trivial, which reflects a
similar situation in general bundles, considered in
Sect.~\ref{Sect7.1}. However, the really interesting and
sensible case is when one considers frames compatible with the
covariant derivatives. As we
know (see~\eref{6.32}), it corresponds to arbitrary non\ndash
degenerate matrix\ndash valued functions $[B_\mu^\nu]$ and
$B=[B_b^a]$ and a matrix\ndash valued functions $B_\mu=[B_{b\mu}^{\nu}]$
given by
    \begin{equation}    \label{7.16}
B_\mu=\tilde{E}_\mu(B)\cdot B^{-1} = B_\mu^\nu E_\nu(B) \cdot B^{-1} .
    \end{equation}
In particular, such are all holonomic frames in $T(E)$, locally
induced by local coordinates on $E$.
Now the normal frames equation~\eref{7.12} (or~\eref{7.1}) reduces to
    \begin{equation}    \label{7.17}
(\Gamma_\mu\cdot B +E_\mu(B))|_{\pi(U)} = 0 .
    \end{equation}
This equation leaves the frame $\{\tilde{E}_\mu= \pi_*|_{\Delta^h}(X_\mu)\}$
in $T(M)$ completely arbitrary and imposes restriction on the frame
$\{\tilde{E}_a=v^{-1}(X_a)=B_a^b E_b\}$ in $E$. This conclusion justifies the
following definition.

    \begin{Defn}    \label{Defn7.3}
Given a \emph{linear} connection $\Delta^h$ on a \emph{vector} bundle
$(E,\pi,M)$ and a subset $U_M\subseteq M$. A \emph{frame $\{E_a\}$ in}
 $E$, defined over an open set $V_M$ containing $U_M$ or equal to it,
$V_M\supseteq U_M$,  is called \emph{normal for $\Delta^h$ over/on}
$U_M$ if their is a frame $\{X_{I}\}$ in $T(E)$, defined over
an open set $V_E\subseteq E$, which is normal for $\Delta^h$ over
a subset $U_E\subseteq E$ and such that $\pi(U_E)=U_M$,
$\pi(V_E)=V_M$, and $E_a=v^{-1}(X_a)$, with the mapping $v$
defined by~\eref{4.1-1}. Respectively, $\{E_a\}$ is normal for
$\Delta^h$ along a mapping $g\colon Q_M\to M$,
$Q_M\not=\varnothing$, if $\{E_a\}$ is normal for $\Delta^h$ over
$g(Q_M)$.
    \end{Defn}

    Taking into account definition~\ref{Defn7.1}, we see that the
so-defined normal frames in the bundle space $E$ are just the ones used in
the theory of frames normal for linear connections in vector
bundles~\cite{bp-NF-D+EP,bp-NF-LTP,bp-Frames-n+point,bp-Frames-path,
bp-Frames-general}.

    It is quite clear, to any frame $\{X_{I}\}$ in $T(E)$ normal over
$U\subseteq E$, there corresponds a unique frame
$\{E_a=v^{-1}(X_a)\}$ in $E$ normal over $\pi(U)\subseteq M$. But, to a
frame  $\{E_a\}$ in $E$ normal over $\pi(U)$, there correspond
infinitely many frames
 $\{X_{I}\}=\{(\pi_*|_{\Delta^h})^{-1}(E_\mu),v(E_a)\}$
in $T(E)$ normal over $U$, where $\{E_\mu\}$ is an arbitrary frame in
$T(M)$ over $\pi(U)$. Thus the problems of existence and
(un)uniqueness of normal frames in $T(E)$ is completely reduced to
the same problems for normal frames in $E$. The last kind of
problems, as we noted at the beginning of the present section,
are known and investigated and the reader is referred
to~\cite{bp-NF-D+EP,bp-NF-LTP,bp-Frames-n+point,bp-Frames-path,
bp-Frames-general} for their solutions and further details.

    Ending, we emphasize that a normal frame $\{E_a\}$ in $E$, as well as
the basis $\{v(E_a)\}$ for $\Delta^v$, can be holonomic as well as
anholonomic (see \emph{loc.~cit.}); at the same time, a normal frame
$\{X_{I}\}$ in $T(E)$ is anholonomic unless some conditions hold%
, a necessary condition being the flatness
(integrability) of the horizontal distribution $\Delta^h$.


\section {Conclusion}
\label{Conclusion}

    In Sect.~\ref{Sect7.1}, we saw that the theory of normal frames in the
most general case is quite trivial. This reflects the understanding that the
more general a concept is, the less particular properties it has, but the more
concrete applications it can find if it is restricted somehow. This situation
was demonstrate when normal frames adapted to  holonomic  ones were considered;
\eg they exist at a given point or along an injective horizontal path, but on
an open set they may exist only in the flat case. A feature of a vector
bundle $(E,\pi,M)$ is that the frames in $T(E)$ over $E$ are in bijective
correspondence with pairs of frames in $E$ over $M$ and in $T(M)$ over $M$.
This result allows the normal frames in $T(E)$, if any, to be `lowered' to
ones in $E$. From here a conclusion was made that the theory of frames in
$T(E)$ normal for linear connections on a vector bundle is equivalent to the
existing one of frames in $E$ normal for covariant derivatives in
$(E,\pi,M)$~\cite{bp-NF-D+EP,bp-NF-LTP}.

    It should be emphasized, the importance of the normal frames for the
physics comes from the fact that they are the mathematical object
corresponding to the physical concept of inertial frame of
reference~\cite{bp-PE-P?,bp-NF-D+EP,bp-EPinED}.


\appendix

\section*
{Appendix: Coordinates normal along injective mappings
					with non-vanishing horizontal component}
\addcontentsline{toc}{section}%
{Appendix: Coordinates normal along injective mappings
					with non-vanishing horizontal component}
\label{Appendix}
\renewcommand{\thesection}{A}

	The purpose of this Appendix is a multi-dimensional generalization of
proposition~\ref{Prop7.5} in the real case, $\field=\field[R]$. It is
formulated below as theorem~\ref{ThmA.1}. For its proof we shall need a
result which is a multidimensional generalization of lemma~\ref{2-Lem6.1-0}.

    \begin{Lem}    \label{3-Lem8.1-0}
Let $n\in\field[N]$, $M$ be a $C^3$ manifold with $\dim M\ge n$, $J^n$ be an
open set in $\field[R]^n$, and
$\gamma \colon J^n\to M$ be $C^1$ regular injective mapping. For every $s_0\in
J^n$, there exists a chart $(U_1,x)$ of $M$ such that $\gamma(s_0)\in U_1$ and
$x(\gamma(s))=(s,\bs{t}_0)$ for some fixed $\bs{t_0}\in\field[R]^{\dimR M-n}$
and all $s\in J^n$ such that $\gamma(s)\in U_1$.
    \end{Lem}

    \begin{Proof}
	Let us choose arbitrary some $s_0\in J^n$ and a chart $(U,y)$ with
$U\ni\gamma(s_0)$ and $y\colon U\to\field[R]^{\dimR M}$. Since the regularity of
$\gamma$ at $s_0$ means that
 $\bigl[\frac{\pd\gamma_y^i}{\pd s^a}\big|_{s_0}\bigr]$
has maximal rank, equal to $n$, we, without loss of generality, can suppose
the coordinates $\{y^i\}$ to be taken such that
 $\det \bigl[\frac{\pd\gamma_y^a}{\pd s^b}\big|_{s_0}\bigr] \not=0,\infty$.%
\footnote{%
If we start from a chart $(U,z)$ for which the matrix
 $ \bigl[{\pd\gamma_z^a}/{\pd s^b}\big|_{s_0}\bigr] $
is degenerate, we can make a coordinate change $\{z^i\}\to\{y^i\}$ with
 $y^i=z^{\alpha_i}$, where the integers $\alpha_1,\dots,\alpha_{\dimR M}$
form a permutation of $1,\dots,\dimR M$, such that
 $ \bigl[{\pd\gamma_y^a}/{\pd s^b}\big|_{s_0}\bigr] $
is non\ndash degenerate. (For the proof, see any book on matrices,
e.g.~\cite{Bellman,Gantmacher/matrices-1}.) Further, we suppose that such a
renumbering of the local coordinates is already done if required. (Cf.\
footnote~\ref{2-renumbering}).%
}
Then the implicit function
theorem~\cite{Warner,Schwartz/Analysis-1,Dieudonne} implies the existence of
a subneighborhood $J_1^n\subseteq J^n$ with $J_1^n\ni s_0$ and such that the
matrix
 $\bigl[ \frac{\pd \gamma_y^a}{\pd s^b}\big|_s \bigr]$
is non\ndash degenerate for $s\in J_1^n$ and the mapping
\[
(\gamma_y^1,\dots,\gamma_y^n)|_{J_1^n}
 \colon J_1^n \to
(\gamma_y^1(J_1^n),\dots,\gamma_y^n(J_1^n)) \subseteq \field[R]^n ,
\]
with
\(
(\gamma_y^1,\dots,\gamma_y^n)|_{J_1^n}\colon s\mapsto
(\gamma_y^1(s),\dots,\gamma_y^n(s))
\)
for $s\in J_1^n$. is a $C^1$ diffeomorphism. Define a chart $(U_1,x)$ of $M$
with domain
	\begin{subequations}		\label{3-8.3}
	\begin{equation}	\label{3-8.3a}
	\begin{split}
U_1 :=& \{ p|p\in U,\ y^a(p)\in\gamma_y^a(J_1^n),\ a=1,\dots,n \}
\\
     =& y^{-1} \bigl(
	(\gamma_y^1(J_1^n),\dots,\gamma_y^n(J_1^n))
	\times \field[R]^{\dimR M-n}
       \bigr)
\ni\gamma(s_0)
	\end{split}
	\end{equation}
and local coordinate functions $x^i$ given via
	\begin{equation}	\label{3-8.3b}
	\begin{split}
y^a &=: \bigl( \gamma_y^a\big|_{J_1^n} \bigr) \circ (x^1,\dots,x^n),
\qquad a=1,\dots,n
\\
y^k &=: x^k +
\bigl( \gamma_y^k\big|_{J_1^n} \bigr) \circ (x^1,\dots,x^n) - t_0^k,
\qquad k=n+1,\dots,\dimR M
	\end{split}
	\end{equation}
	\end{subequations}
where $(x^1,\dots,x^n)\colon p\mapsto (x^1(p),\dots,x^n(p))$, $p\in U_1$, and
$t_0^k\in\field[R]$ are constant numbers.  Since
 $\frac{\pd y^a}{\pd x^b} = \frac{\pd\gamma_y^a}{\pd s^b}$,
 $\frac{\pd y^a}{\pd x^k} = \delta_k^a$ for $k\ge n+1$,
 $\frac{\pd y^k}{\pd x^a} = \frac{\pd\gamma_y^k}{\pd s^a}$ for $k\ge n+1$,
and
 $\frac{\pd y^k}{\pd x^l} = \delta_l^k$ for $k,l\ge n+1$,
the Jacobian of the change $\{y^i\}\to\{x^i\}$ on $U_1$ is
\(
\det\bigl[\frac{\pd x^i}{\pd y^j}\bigr]
=
\bigl( \det\bigl[\frac{\pd \gamma_y^a}{\pd s^b}\bigr] \bigr)^{-1}
\not= 0,\infty .
\)
Consequently $x^i$ are really coordinate functions and
$x\colon U_1\to J_1^n\times\field[R]^{\dimR M-n}$
is in fact coordinate homeomorphism.%
\footnote{%
The so-constructed chart $(U_1,x)$ is, obviously, a multidimensional
generalization of a similar chart defined in the proof of
lemma~\ref{2-Lem6.1-0} --- see  the paragraph containing
equation~\eref{2-6.20}.%
}
The coordinates $\{x^i\}$ can be expressed through $\{y^i\}$ explicitly.
Indeed, writing the first raw of~\eref{3-8.3b} as
	\begin{multline*}
(y^1,\dots,y^n)
= (\gamma_y^1|_{J_1^n},\dots,\gamma_y^n|_{J_1^n})\circ (x^1,\dots,x^n)
= (\gamma_y^1,\dots,\gamma_y^n)|_{J_1^n}\circ (x^1,\dots,x^n)
	\end{multline*}
and using that $(\gamma_y^1,\dots,\gamma_y^n)|_{J_1^n}$ is a $C^1$
diffeomorphism and the second raw of~\eref{3-8.3b}, we find
(cf.~\eref{2-6.20})
	\begin{equation}
	\tag{\ref{3-8.3b}$^\prime$}	\label{3-8.3b'}
	\begin{split}
&
(x^1,\dots,x^n)
= \bigl( (\gamma_y^1,\dots,\gamma_y^n)|_{J_1^n} \bigr)^{-1}
  \circ (y^1,\dots,y^n)
\\
&
x^k = y^k
- (\gamma_y^k|_{J_1^n})
  \circ \bigl( (\gamma_y^1,\dots,\gamma_y^n)|_{J_1^n} \bigr)^{-1}
  \circ (y^1,\dots,y^n)
+ t_0^k,
\qquad k\ge n+1  .
	\end{split}
	\end{equation}

	Using~\eref{3-8.3}, we see that in $(U_1,x)$ the local coordinates
of $\gamma(s)$ for $s=(s^1,\dots,s^n)\in J_1^n$ are
	\begin{equation}	\label{3-8.4}
\gamma^a(s):=x^a(\gamma(s)) = s^a,
\quad
\gamma^k(s):=x^k(\gamma(s)) = t_0^k,
\qquad k\ge n+1,
	\end{equation}
\ie $x(\gamma(s))=(s,\bs{t}_0)$ for some fixed
$\bs{t}_0=(t_0^{n+1},\ldots,t_0^{\dimR M})\in\field[R]^{\dimR M-n}$.
    \end{Proof}

	Thus, in the chart $(U_1,x)$ or the coordinates $\{x^i\}$ constructed
above, the first $n$ coordinates of a point lying in $\gamma(J^n)$, \ie in
$\gamma(J_1^n)$, coincide with the corresponding parameters $s^1,\ldots,s^n$ of
$\gamma$, the remaining coordinates, if any, being constant numbers. This
conclusion allows locally, in $U_1$, the mapping $\gamma$ to be considered as a
representative of a family of mappings
$\eta(\cdot,\bs{t})\colon J_1^n\to M$, $\bs{t}\in\field[R]^{\dimR M-n}$,
defined by $\eta(s,\bs{t}):=x^{-1}(s,\bs{t})$ for
 $(s,\bs{t})\in J_1^n\times\field[R]^{\dimR M-n}$.
In fact, we have
$\gamma=\eta(\cdot,\bs{t}_0)$ or $\gamma(s)=\eta(s,\bs{t}_0)$.%
\footnote{%
In~\cite{bp-Frames-general} the existence of $\eta$ is taken as a given fact
without proof.%
}

	Let $(E,\pi,M)$ be a $C^3$ bundle endowed with  $C^1$ connection
$\Delta^h$. Let $k\in\field[N]$, $k\le\dim M$, and $J^k$ be an open set in
$\field[R]^k$.
	Consider a $C^2$ regular injective mapping $\beta\colon J^k\to E$ such
that the vector fields
\(
\dot{\beta}_\alpha\colon s \mapsto
\dot{\beta}_\alpha(s) := \frac{\pd\beta^I(s)}{\pd s^\alpha}
				\frac{\pd}{\pd u^I}\big|_{\beta(s)},
\)
with $s:=(s^1,\dots,s^k)\in J^k$ and $\alpha=1,\dots,k$,
do \emph{not} belong to the vertical distribution $\Delta^v$,
$\dot{\beta}_\alpha(s) \not \in\Delta_{\beta(s)}^v$ for all $s\in J^k$; in
particular, the mapping $\beta$ can be a horizontal mapping in a sense that
$\dot{\beta}_\alpha(s) \in\Delta_{\beta(s)}^h$ for all $s\in J^k$, but
generally these vectors can have a vertical component too.
	Our aim is to find the integrability conditions for the normal
frame/coordinates equation~\eref{7.6} and its solutions, if any, when
$U=\beta(J^k_1)$ for some subset $J^k_1\subseteq J^k$.

	Let us take some $s_0\in J^k$ and construct the chart $(U_1,u)$ with
$U_1\ni\beta(s_0)$ provided by lemma~\ref{3-Lem8.1-0} with $E$ for $M$ and
$\beta$ for $\gamma$.
	If $J_1^k:=\{s\in J^k : \beta(s)\in U_1\}$ and $p\in U_1$, then there is
a unique $(s,\bs{t})\in J_1^k\times\field[R]^{\dimR E-k}$ such that
$p=\eta(s,\bs{t})$ with $\eta:=u^{-1}$, \ie
$u^I(p)=s^I$ for $I=1,\dots,k$ and
$u^I(p)=t^I$ for $I=k+1,\dots,n+r$.
Besides, we have $u(\beta(s))=(s,\bs{t}_0)$ for all $s\in J_1^k$ and some fixed
$\bs{t}_0\in\field[R]^{\dimR E-k}$. Since the vector fields
$\dot{\beta}_\alpha$, $\alpha=1,\dots,k$, are not vertical, we can construct
the coordinates $\{u^I\}$, associated to the chart $(U_1,u)$, so that they to
be \emph{bundle} coordinates on $U_1$ (see the proof of lemma~\ref{3-Lem8.1-0}).
	Thus on $U_1$ we have bundle coordinates $\{u^I\}$ such that
	\begin{equation}	\label{A.1}
	\begin{split}
\big( u^1(\eta(s,\bs{t})),\dots,u^{n+r}(\eta(s,\bs{t})) \big)
:= (s,\bs{t})\in\field[R]^{n+r}
 \\
s=(s^1,\dots,s^k)\in J_1^k
\quad
\bs{t}=(t^{k+1},\dots,t^{n+r}) \in \field[R]^{n+r-k} .
	\end{split}
	\end{equation}

	Let the indices $\alpha$ and $\beta$ run from $1$ to $k$ and the
indices $\sigma$ and $\tau$ take the values form $k+1$ to $n$; we set
$\sigma=\tau=\varnothing$ if $k=n$. Thus, we have
$u^\alpha(\eta(s,\bs{t}))=s^\alpha$,
$u^\sigma(\eta(s,\bs{t}))=t^\sigma$, and
$u^a(\eta(s,\bs{t}))=t^a$.

	\begin{Prop}	\label{PropA.1}
Under the hypotheses made above, the normal frame/coordinates
equation~\eref{7.6} with $U=\beta(J_1^k)=\beta(J^k)\cap U_1$ has solutions if
and only if the system of equations
	\begin{equation}	\label{A.2}
\Bigl(
  \frac{\pd\Gamma_\alpha^b}{\pd u^\beta}
- \frac{\pd\Gamma_\beta^b}{\pd u^\alpha}
\Bigr)\Big|_{\beta(s)} B_b^a(s)
+
\Gamma_\alpha^b(\beta(s)) \frac{\pd B_b^a(s)}{\pd s^\beta} -
\Gamma_\beta^b(\beta(s)) \frac{\pd B_b^a(s)}{\pd s^\alpha}
= 0 ,
	\end{equation}
where $\Gamma_\mu^a$ are the 2-index coefficients of $\Delta^h$ in $\{u^I\}$,
has solutions $B_b^a\colon J_1^k\to\field[R]$ with $\det[B_b^a]\not=0,\infty$.
Besides, if such solutions exist, then all solutions of~\eref{7.6}  are given
on $U_1$ by the formula
	\begin{multline}	\label{A.3}
\tilde{u}^a(\eta(s,\bs{t}))
=
- \int\limits_{s_1}^{s} B_b^a(s) \Gamma_\alpha^b(\beta(s)) \Id s^\alpha
\\
-
B_b^a(s) \Gamma_\mu^b(\beta(s)) [u^\mu(\eta(s,\bs{t})) - u^\mu(\beta(s))]
+
B_b^a(s)[u^b(\eta(s,\bs{t})) - u^b(\beta(s))]
\\ +
f_{\mu\nu}^a(s;\bs{t};\eta)
[u^\mu(\eta(s,\bs{t})) - u^\mu(\beta(s))]
[u^\nu(\eta(s,\bs{t})) - u^\nu(\beta(s))] ,
	\end{multline}
where $s_1\in J_1^k$ is arbitrarily fixed, $B_b^a$, with
$\det[B_b^a]\not=0,\infty$, are solutions of~\eref{A.2}, and the functions
$f_{\mu\nu}^a$ and their first partial derivatives are bounded when
$\bs{t}\to\bs{t}_0$.
	\end{Prop}

	\begin{Rem}	\label{RemA.1}
As $u^\alpha(\eta(s,\bs{t}))=u^\alpha(\beta(s))\equiv s^\alpha$ for all
$\alpha=1,\dots,k$, the terms with $\mu,\nu=1,\dots,k$ in~\eref{A.3} have
vanishing contribution.
	\end{Rem}

	\begin{Rem}	\label{RemA.2}
For $k=1$, we have $\alpha=\beta=1$, due to which the equations~\eref{A.2}
are identically valid and proposition~\ref{PropA.1} reduces to
proposition~\ref{Prop7.5}.
	\end{Rem}

	\begin{Proof}
To begin with, we rewrite~\eref{7.6} as
\[
\frac{\pd \tilde{u}^a}{\pd s^\alpha} \Big|_{\beta(s)}
=
- \frac{\pd \tilde{u}^a}{\pd t^b} \Big|_{\beta(s)} \Gamma_\alpha^b(\beta(s))
\quad
\frac{\pd \tilde{u}^a}{\pd t^\sigma} \Big|_{\beta(s)}
=
- \frac{\pd \tilde{u}^a}{\pd t^b} \Big|_{\beta(s)} \Gamma_\sigma^b(\beta(s)) .
\]
Introducing a non-degenerate matrix-valued function $[B_b^a]$ on $J_1^k$ by
	\begin{equation}	\label{A.4}
B_b^a(s)
= \frac{\pd \tilde{u}^a}{\pd t^b} \Big|_{\beta(s)}
= \frac{\pd \tilde{u}^a(s,\bs{t})}{\pd t^b} \Big|_{\bs{t}=\bs{t}_0},
	\end{equation}
we see that~\eref{7.6} is equivalent to
	\begin{subequations}	\label{A.5}
	\begin{align}	\label{A.5a}
\frac{\pd \tilde{u}^a}{\pd s^\alpha} \Big|_{\beta(s)}
& =
- B_b^a(s) \Gamma_\alpha^b(\beta(s)) \quad  \alpha=1,\dots,k
\\			\label{A.5b}
\frac{\pd \tilde{u}^a}{\pd t^\sigma} \Big|_{\beta(s)}
& =
- B_b^a(s) \Gamma_\sigma^b(\beta(s)) \quad \sigma=k+1,\dots,n .
	\end{align}
	\end{subequations}

	Expanding $\tilde{u}^a(\eta(s,\bs{t}))$ into a Taylor's polynomial up
to second order terms relative to $(\bs{t}-\bs{t}_0)$ about the point
$\bs{t}_0$ and using~\eref{A.4} and~\eref{A.5}, we get:
	\begin{multline}	\label{A.6}
\tilde{u}^a(\eta(s,\bs{t}))
=
f^a(s)
- B_b^a(s) \Gamma_\sigma^b(\beta(s)) [t^\sigma-t_0^\sigma]
+ B_b^a(s)[t^b-t_0^b]
+ f_{\sigma\tau}^a(s;\bs{t};\eta) [t^\sigma-t_0^\sigma] [t^\tau-t_0^\tau]
\\
=
f^a(s)
-
B_b^a(s) \Gamma_\mu^b(\beta(s)) [u^\mu(\eta(s,\bs{t})) - u^\mu(\beta(s))]
+
B_b^a(s)[u^b(\eta(s,\bs{t})) - u^b(\beta(s))]
\\ +
f_{\mu\nu}^a(s;\bs{t};\eta)
[u^\mu(\eta(s,\bs{t})) - u^\mu(\beta(s))]
[u^\nu(\eta(s,\bs{t})) - u^\nu(\beta(s))] ,
	\end{multline}
where $f^a$ and $f_{\mu\nu}^a$ are $C^1$ functions and $f_{\mu\nu}^a$ and
their first partial derivatives are bounded when $\bs{t}\to\bs{t}_0$. The
equation~\eref{A.5a} is the only condition that puts some restrictions on
$f^a$ and $B_b^a$ (besides $\det[B_b^a]\not=0,\infty$). Inserting~\eref{A.6}
into~\eref{A.5a} and using that $\beta(s)=\eta(s,\bs{t}_0)$, we obtain
	\begin{equation}	\label{A.7}
\frac{\pd f^a(s)}{\pd s^\alpha}
= - B_b^a(s) \Gamma_\alpha^b(\beta(s)) .
	\end{equation}
Thus the initial normal coordinates equation~\eref{7.6}, with $U=\beta(J_1^k)$,
has solutions if and only if there exist solutions of~\eref{A.7} relative to
$f^a$ and/or $B_b^a$. The integrability conditions for~\eref{A.7}
are~\cite{Hartman}
\[
0 =
\frac{\pd^2 f^a}{\pd s^\beta \pd s^\alpha} -
\frac{\pd^2 f^a}{\pd s^\alpha \pd s^\beta}
=
- \frac{\pd}{\pd s^\beta} \bigl( B_b^A(s) \Gamma_\alpha^b(s)) \bigr)
+ \frac{\pd}{\pd s^\alpha} \bigl( B_b^A(s) \Gamma_\beta^b(s)) \bigr)
= \dotsb
\]
and coincide with~\eref{A.2}, by virtue of $u^\alpha(\beta(s))=s^\alpha$.
This result concludes the proof of the fires part of the proposition.

	If~\eref{A.2} admits solutions $B_b^a$ with
$\det[B_b^a]\not=0,\infty$, then the general solution of~\eref{A.7} is
\(
f^a(s)
= - \int\limits_{s_1}^{s} B_b^a(s) \Gamma_\alpha^b(\beta(s)) \Id s^\alpha
\)
for some $s_1\in J_1^k$ and this solution is independent of the integration
path in $J_1^k$, due to~\eref{A.2}.
	\end{Proof}

	\begin{Lem}	\label{LemA.1}
Let $(E,\pi,M)$ be a $C^3$ bundle endowed with  $C^2$ connection with
coefficients $\Gamma_\mu^a$ in the frame adapted to local coordinates
$\{u^i\}$, defined before proposition~\ref{PropA.1}. There exist solutions
$B_b^a$ with $\det[B_b^a]\not=0,\infty$ of the system of equations~\eref{A.2}
if and only if the coefficients $\Gamma_\mu^a$ satisfy the equations
	\begin{subequations}	\label{A.12}
	\begin{align}	\label{A.12a}
R_{\alpha\beta}^{a} (\beta(s)) & = 0 \qquad s\in J_1^k
\\			\label{A.12b}
\Bigl(
    \Gamma_\alpha^d \frac{\pd^2\Gamma_\beta^c}{\pd u^b\pd u^d}
  - \Gamma_\beta^d \frac{\pd^2\Gamma_\alpha^c}{\pd u^b\pd u^d}
\Bigr) \Big|_{\beta(s)}
& = 0
\qquad s\in J_1^k
	\end{align}
	\end{subequations}
in which $R_{\mu\nu}^{a}$ are the (fibre) components in $\{u^I\}$ of the
curvature of $\Delta^h$, defined by~\eref{3.24a}. If the
conditions~\eref{A.12} are valid, the set of the solutions of~\eref{A.2}
coincides with the set of solutions of the system
	\begin{equation}	\label{A.13}
\frac{\pd B_b^a(s)}{\pd s^\alpha}
=
- B_c^a(s) \frac{\pd \Gamma_\alpha^c}{\pd u^b} \Big|_{\beta(s)}
+ \frac{\pd D_{b}^{a}(s)}{\pd s^\alpha}
	\end{equation}
relative to $B_b^a$, where $D_b^a$ are solutions of
	\begin{equation}	\label{A.14}
\Bigl(
 \Gamma_\alpha^b(\beta(s)) \frac{\pd}{\pd s^\beta}
- \Gamma_\beta^b(\beta(s)) \frac{\pd}{\pd s^\alpha}
\Bigr) D_b^a(s)
= 0 .
	\end{equation}
	\end{Lem}

	\begin{Proof}
	Consider the integrability condition~\eref{A.2} for~\eref{7.6} in
more details. Define functions $D_{b\alpha}^a\colon J_1^k\to\field=\field[R]$
via the equation
	\begin{equation}	\label{A.8}
\frac{\pd B_b^a(s)}{\pd s^\alpha}
=
- B_c^a(s) \frac{\pd \Gamma_\alpha^c}{\pd u^b} \Big|_{\beta(s)}
+ D_{b\alpha}^{a}(s) .
	\end{equation}
The substitution of this equality into~\eref{A.2} results in
\[
R_{\beta\alpha}^{b}(\beta(s)) B_b^a(s)
- \Gamma_\alpha^b(\beta(s)) D_{b\beta}^{a}(s)
+ \Gamma_\beta^b(\beta(s)) D_{b\alpha}^{a}(s)
= 0 ,
\]
where $R_{\alpha\beta}^{a}$ are the (fibre) components in $\{u^I\}$ of the
curvature of $\Delta^h$, defined by~\eref{3.24a}.
The simple observation that $\{\tilde{u}^\alpha,\tilde{u}^a\}$, if they
exist as solutions of~\eref{7.6}, are normal coordinates on the whole bundle
space of the restricted bundle $(U,\pi|_U,\pi(U))$ with $U=\beta(J_1^k)$ leads
to
	\begin{equation}	\label{A.9}
R_{\alpha\beta}^{a} (\beta(s)) = 0 \qquad s\in J_1^k  ,
	\end{equation}
by virtue of proposition~\ref{Prop7.6}. Therefore the previous equation
reduces to
	\begin{subequations}	\label{A.10}
	\begin{equation}	\label{A.10a}
  \Gamma_\alpha^b(\beta(s)) D_{b\beta}^{a}(s)
- \Gamma_\beta^b(\beta(s)) D_{b\alpha}^{a}(s)
= 0 .
	\end{equation}
It is clear that~\eref{A.8}--\eref{A.10a} are equivalent to~\eref{A.2}.
Consequently, the quantities $D_{b\alpha}^{a}$ must be solutions
of~\eref{A.10a} while the $C^1$ functions $B_b^a$ have to be solutions
of~\eref{A.8}. The integrability conditions
\(
\bigl( \frac{\pd^2}{\pd s^\beta \pd s^\alpha}
     - \frac{\pd^2}{\pd s^\alpha \pd s^\beta}
\bigr) B_b^a(s) = 0
\)
for~\eref{A.8} can be written as~%
\footnote{~%
At this point one should require $\Delta^h$ to be of class $C^2$ which is
possible if the manifolds $E$ and $M$ are of class $C^3$.%
}
	\begin{multline*}
\Bigl(
- \frac{\pd^2\Gamma_\alpha^c}{\pd u^\beta \pd u^b}
+ \frac{\pd^2\Gamma_\beta^c}{\pd u^\alpha \pd u^b}
+ \frac{\pd\Gamma_\alpha^d}{\pd u^b} \frac{\pd\Gamma_\beta^c}{\pd u^d}
- \frac{\pd\Gamma_\beta^d}{\pd u^b} \frac{\pd\Gamma_\alpha^c}{\pd u^d}
\Bigr) \Big|_{\beta(s)}
B_c^a(s)
+
  \frac{\pd D_{b\alpha}^{a}(s)}{\pd s^\beta}
- \frac{\pd D_{b\beta}^{a}(s)}{\pd s^\alpha}
= 0
	\end{multline*}
which conditions split into
	\begin{equation}	\label{A.10b}
0 =
  \frac{\pd D_{b\alpha}^{a}(s)}{\pd s^\beta}
- \frac{\pd D_{b\beta}^{a}(s)}{\pd s^\alpha}
	\end{equation}
	\end{subequations}
\vspace{-4ex}
	\begin{multline}	\label{A.11}
 0 =
\Bigl(
- \frac{\pd^2\Gamma_\alpha^c}{\pd u^\beta \pd u^b}
+ \frac{\pd^2\Gamma_\beta^c}{\pd u^\alpha \pd u^b}
+ \frac{\pd\Gamma_\alpha^d}{\pd u^b} \frac{\pd\Gamma_\beta^c}{\pd u^d}
- \frac{\pd\Gamma_\beta^d}{\pd u^b} \frac{\pd\Gamma_\alpha^c}{\pd u^d}
\Bigr) \Big|_{\beta(s)}
=
\Bigl(
  - \Gamma_\alpha^d \frac{\pd^2\Gamma_\beta^c}{\pd u^b\pd u^d}
  + \Gamma_\beta^d \frac{\pd^2\Gamma_\alpha^c}{\pd u^b\pd u^d}
\Bigr) \Big|_{\beta(s)} ,
	\end{multline}
where~\eref{A.9} and~\eref{3.24a} were applied in the derivation of the
second equality in~\eref{A.11}.

	Since the system of equations~\eref{A.10} always has solutions, \eg
$D_{b\alpha}^{b}(s)=0$, we can assert that~\eref{A.9} and~\eref{A.11} are the
integrability conditions for~\eref{A.2} and, if~\eref{A.9} and~\eref{A.11}
hold, every solution of~\eref{A.8}, with $D_{b\alpha}^{a}$
satisfying~\eref{A.10}, is a solution of~\eref{A.2} and \emph{vice versa}.

	At the end, the only unproved assertion is that $D_{b\alpha}^{a}$
in~\eref{A.8} equals to $\pd_\alpha(D_b^a)$ with $D_b^a$ satisfying~\eref{A.14}.
Indeed, since $J_1^k$ is an open set and hence is contractible one, the
Poincar\'e's lemma (see~\cite[sec.~6.3]{Nash&Sen}
or~\cite[pp.~21,~106]{Gockeler&Schucker}) implies the existence of functions
$D_b^a$ on $J_1^k$ such that $D_{b\alpha}^{b}(s)=\pd_\alpha(D_b^a)(s)$, due
to~\eref{A.10b}; inserting this result into~\eref{A.10a}, we get~\eref{A.14}.
	\end{Proof}

	\begin{Rem}	\label{RemA.3}
Regardless that the conditions~\eref{A.12b} look quite special, they are
identically valid for connections with
	\begin{equation}	\label{A.15}
\Gamma_\alpha^a
= - (\Gamma_{b\alpha}^{a}\circ\pi)\cdot u^b + G_\alpha^a\circ\pi ,
	\end{equation}
where $\Gamma_{b\alpha}^{a}$ and $G_\alpha^a$ are $C^2$ functions on
$\pi(\beta(J^k))$. In particular, of this kind are the linear connections on
vector bundles --- see proposition~\ref{Prop6.1}.
	\end{Rem}

	At last, we shall formulate the main result of the above
considerations as a combination of proposition~\ref{PropA.1} and
lemma~\ref{LemA.1}.

	\begin{Thm}	\label{ThmA.1}
Let $(E,\pi,M)$ be a $C^3$ bundle endowed with a $C^2$ connection. Under the
hypotheses made and notation introduced before proposition~\ref{PropA.1},
there exist solutions of the normal frame/coordinates equation~\eref{7.6} if
and only if the connection's coefficients satisfy the equations~\eref{A.12}.
If these equations hold, all coordinates normal on $\beta(J_1^k)$ are given on
$U_1$ by~\eref{A.3}, where $B_b^a$ are solutions of~\eref{A.13}, with $D_b^a$
being solutions of~\eref{A.14}.
	\end{Thm}

    \begin{Rem}    \label{RemA.4}
If there are $s_0\in J^k$ and $\alpha\in\{1,\dots,k\}$ such that the vector
$\dot{\beta}_\alpha(s_0)$ is a vertical vector,
$\dot{\beta}_\alpha(s_0)\in\Delta_{\beta(s_0)}^v$, then theorem~\ref{ThmA.1}
remains true with the only correction that the coordinates $\{u^I\}$ will
\emph{not be bundle} coordinates. If this is the case, the constructed
coordinates $\{\tilde{u}^I\}$ will be solutions of~\eref{7.6}, but we cannot
assert that they are bundle coordinates which are (locally) normal along
$\beta$ in a neighborhood of the point $\beta(s_0)$.
    \end{Rem}

	Theorem~\ref{ThmA.1} provides a necessary and sufficient condition for
the existence of local coordinates in a neighborhood of $\beta(s_0)$ for any
$s_0\in J^k$ which are locally normal along $\beta$, \ie on $\beta(J_1^k)$ for
some open subset $J_1^k\subseteq J^k$ containing $s_0$. Moreover, if this
condition is valid, the theorem describes locally all coordinates normal along
$\beta$.

	Theorem~\ref{ThmA.1} can be generalized by requiring $\beta$ to be
locally injective instead of injective, \ie for each $s\in J$ to exist subset
$J_s^k\subseteq J^k$ such that $J_s^k\ni s$ and the restricted
mapping $\beta|_{J_s^k}$ to be injective.
	Besides, if one needs a version of the above results for complex
bundles, they should be considered as real ones (with doubled dimension of the
manifolds) for which are applicable the above considerations.




\addcontentsline{toc}{section}{References}
\bibliography{bozhopub,bozhoref}

\begin{thebibliography}{10}

\bibitem{Riemann/Hypotheses}
B.~Riemann.
\newblock {\"{U}ber die Hypothesen welche de Geometrie zugunde liegen} ({On}
  the hypotheses underlying the geometry).
\newblock {\em {G{\"o}tingen} Abh.}, 13:1--20, 1868.
\newblock Habilitationsschrifr (Ph.d.\ thesis), 1854. In German.

\bibitem{Veblen}
O.~Veblen.
\newblock Normal {Coordinates} for the {Geometry} of {Paths}.
\newblock {\em Proc.\ {N}at.\ {A}cad.}, 8:192--197, 1922.

\bibitem{Fermi}
E.~Fermi.
\newblock Sopra i fenomeni che avvengono in vicinonza di una linear oraria
  ({On} phenomena near a world line).
\newblock {\em Atti {R}.\ {A}ccad {L}incei {R}end., {C}l.\ {S}ci.\ {F}is.\
  {M}at.\ {N}at.}, 31(1):21--23, 51--52, 1922.
\newblock Russian translation: Enrico Fermi, Scientific papers, vol.~I,
  1921--1938. Italy, Nauka, Moscow, 1971, pp.64--71.

\bibitem{Levi-Civita/1926}
T.~Levi-Civita.
\newblock Sur l'\'ecart g\'eod\'esique.
\newblock {\em Math.\ Ann.}, 97:291--320, 1926.

\bibitem{Eisenhart/Non-Riemannian}
Luther~P. Eisenhart.
\newblock {\em Non-{Riemannian} geometry}, volume VIII of {\em Colloquium
  Publications}.
\newblock American Mathematical Society, New York, 1927.

\bibitem{ORai}
L.~\'O~Raifeartaigh.
\newblock Fermi coordinates.
\newblock {\em Proceedings of the {Royal} {Irish} {Academy}}, 59
  {Sec.}~A(2):15--24, 1958.

\bibitem{bp-Frames-n+point}
Bozhidar~Z. Iliev.
\newblock Normal frames and the validity of the equivalence principle: {I}.
  {Cases} in a neighborhood and at a point.
\newblock {\em Journal of Physics A: Mathematical and General},
  29(21):6895--6901, 1996.
\newblock \\ http://www.arXiv.org e-Print archive, E-print No.\ gr-qc/9608019,
  August 1998.

\bibitem{bp-Frames-path}
Bozhidar~Z. Iliev.
\newblock Normal frames and the validity of the equivalence principle: {II}.
  {The} case along paths.
\newblock {\em Journal of Physics A: Mathematical and General},
  30(12):4327--4336, 1997.
\newblock \\ http://www.arXiv.org e-Print archive, E-print No.\ gr-qc/9709053,
  September 1997.

\bibitem{bp-Frames-general}
Bozhidar~Z. Iliev.
\newblock Normal frames and the validity of the equivalence principle: {III}.
  {The} case along smooth maps with separable points of self-intersection.
\newblock {\em Journal of Physics A: Mathematical and General},
  31(4):1287--1296, January 1998.
\newblock \\ http://www.arXiv.org e-Print archive, E-print No.\ gr-qc/9805088,
  May 1998.

\bibitem{bp-NF-D+EP}
Bozhidar~Z. Iliev.
\newblock Normal frames for derivations and linear connections and the
  equivalence principle.
\newblock {\em Journal of Geometry and Physics}, 45(1--2):24--53, February
  2003.
\newblock \\ http://www.arXiv.org e-Print archive, E-print No.\ hep-th/0110194,
  October 2001.

\bibitem{bp-NF-LTP}
Bozhidar~Z. Iliev.
\newblock Normal frames and linear transports along paths in vector bundles.
\newblock \\ http://www.arXiv.org e-Print archive, E-print No.\ gr-qc/9809084,
  September 1998, 2002.

\bibitem{K&N}
S.~Kobayashi and K.~Nomizu.
\newblock {\em Foundations of Differential Geometry}, volume I and II.
\newblock Interscience Publishers, New York-London-Sydney, 1963 and 1969.
\newblock Russian translation: Nauka, Moscow, 1981.

\bibitem{Warner}
F.~W. Warner.
\newblock {\em Foundations of differentiable manifolds and Lie groups}.
\newblock Springer-Verlag, New York-Berlin-Heidelberg-Tokyo, 1983.
\newblock Russian translation: Mir, Moscow, 1987.

\bibitem{Poor}
Walter~A. Poor.
\newblock {\em Differential geometric structures}.
\newblock McGraw-Hill Book Company Inc., New York, 1981.

\bibitem{Spivak-1}
M.~Spivak.
\newblock {\em A comprehensive introduction to differential geometry},
  volume~1.
\newblock Publish or Perish, Boston, 1970.

\bibitem{Mangiarotti&Sardanashvily}
L.~Mangiarotti and G.~Sardanashvily.
\newblock {\em Connections in classical and quantum field theory}.
\newblock World Scientific, Singapore-New Jersey-London-Hong Kong, 2000.

\bibitem{Schutz}
Bernard~F. Schutz.
\newblock {\em Geometrical methods of mathematical physics}.
\newblock Cambridge University Press, Cambridge-London-New York- New
  Rochelle-Melbourne-Sydney, 1982.
\newblock Russian translation: Mir, Moscow, 1984.

\bibitem{Rahula}
Maido Rahula.
\newblock {\em New problems in differential geometry}, volume~8 of {\em Series
  on Soviet and East European Mathematics}.
\newblock World Scientific, Singapore-New Jersey-London-Hong Kong, 1993.

\bibitem{K&N-1}
S.~Kobayashi and K.~Nomizu.
\newblock {\em Foundations of Differential Geometry}, volume~I.
\newblock Interscience Publishers, New York-London, 1963.
\newblock Russian translation: Nauka, Moscow, 1981.

\bibitem{Schwartz/Analysis-1}
Laurent Schwartz.
\newblock {\em Analyse math\'ematique}, volume~I.
\newblock Hermann, Paris, 1967.
\newblock In French; Russian translation: Mir, Moscow, 1972.

\bibitem{Dieudonne}
J.~Dieudonn\'e.
\newblock {\em Foundations of modern analysis}.
\newblock Academic Press, New York, 1960.

\bibitem{bp-PE-P?}
Bozhidar~Z. Iliev.
\newblock Is the principle of equivalence a principle?
\newblock {\em Journal of Geometry and Physics}, 24(3):209--222, 1998.
\newblock \\ http://www.arXiv.org e-Print archive, E-print No.\ gr-qc/9806062,
  June 1998.

\bibitem{bp-EPinED}
Bozhidar~Z. Iliev.
\newblock Equivalence principle in classical electrodynamics.
\newblock In Andrew Chubykalo, Vladimir Onoochin, Augusto Espinoza, and Roman
  Smirnov-Rueda, editors, {\em Has the last word been said on classical
  electrodynamics? -- New Horizons}, pages ??--?? Rinton Press, ??, 2004.
\newblock To appear. \\ http://www.arXiv.org e-Print archive, E-print No.\
  gr-qc/0303002, March 2003.

\bibitem{Bellman}
R.~Bellman.
\newblock {\em Introduction to matrix analysis}.
\newblock McGraw-Hill book comp., New York-Toronto-London, 1960.
\newblock Russian translation: Nauka, Moscow, 1978.

\bibitem{Gantmacher/matrices-1}
F.~R. Gantmacher.
\newblock {\em The theory of matrices}, volume one.
\newblock Chelsea Pub.\ Co., New York, N.Y., 1960.
\newblock Translation from Russian. The right English transliteration of the
  author's name is Gantmakher.

\bibitem{Hartman}
Ph. Hartman.
\newblock {\em Ordinary Differential Equations}.
\newblock John Wiley {\&} Sons, New York-London-Sydney, 1964.

\bibitem{Nash&Sen}
C.~Nash and S.~Sen.
\newblock {\em Topology and Geometry for physicists}.
\newblock Academic Press, London-New York, 1983.

\bibitem{Gockeler&Schucker}
M.~{G{\"o}ckeler} and T.~{Sch\"ucker}.
\newblock {\em Differential geometry, gauge theories, and gravity}.
\newblock Cambridge Univ.\ Press, Cambridge, 1987.

\end{thebibliography}
\bibliographystyle{unsrt}
\addcontentsline{toc}{subsubsection}{This article ends at page}

\end{document}